\documentclass[11pt]{amsart}

\usepackage{a4wide}
\usepackage[margin=1in]{geometry}
\usepackage[skip=4pt plus 1pt minus 1pt]{parskip}
\usepackage[T1]{fontenc}
\usepackage{lmodern}
\usepackage{amsmath,amssymb,amsthm,mathtools}
\usepackage{tikz,booktabs,array,longtable,xcolor}
\usetikzlibrary{arrows.meta}
\usepackage[hidelinks]{hyperref}

\makeatletter
\renewcommand{\l@section}
  {\@tocline{1}{0.7em}{0em}{2.7em}{\bfseries}}

\renewcommand{\l@subsection}
  {\@tocline{2}{0.15em}{1.7em}{3.2em}{}}

\renewcommand{\l@subsubsection}
  {\@tocline{3}{0.1em}{3.5em}{4.0em}{\small}}
\makeatother

\newcommand{\ppmj}{\mathcal{P}_\pi^m(j)}

\newcommand{\E}{\mathbb E}
\newcommand{\Prob}{\mathbb P}
\newcommand{\R}{\mathbb R}

\newcommand{\ind}{\mathbf 1}
\newcommand{\Law}{\operatorname{Law}}

\newtheorem{theorem}{Theorem}[section]
\newtheorem{proposition}[theorem]{Proposition}
\newtheorem{lemma}[theorem]{Lemma}
\newtheorem{corollary}[theorem]{Corollary}

\theoremstyle{definition}
\newtheorem{definition}[theorem]{Definition}
\newtheorem{conjecture}[theorem]{Conjecture}

\theoremstyle{remark}
\newtheorem{remark}[theorem]{Remark}

\numberwithin{equation}{section}

\newcommand{\PackingTreeFigure}{%
\begin{tikzpicture}[x=1.6cm,y=0.8cm,
  vertex/.style={circle,draw,fill=white,inner sep=2.4pt,font=\small},
  trunk/.style={vertex,very thick,fill=black!8}]
\node[trunk] (p1) at (0,1) {$1$};
\node[trunk] (p4) at (0,3) {$4$};
\node[trunk] (p7) at (0,5) {$7$};
\node[trunk] (p10) at (0,7) {$10$};
\node[vertex] (p3) at (-1,0) {$3$};
\node[vertex] (p6) at (-1,2) {$6$};
\node[vertex] (p9) at (-1,4) {$9$};
\node[vertex] (p2) at (1,4) {$2$};
\node[vertex] (p5) at (1,6) {$5$};
\node[vertex] (p8) at (1,8) {$8$};
\draw[very thick] (p1)--(p4)--(p7)--(p10);
\draw (p3)--(p1) (p6)--(p4) (p9)--(p7);
\draw (p4)--(p2) (p7)--(p5) (p10)--(p8);
\draw[-{Stealth},thin] (2.2,2) -- (2.2,5.7)
  node[midway,right,align=left,font=\small] {increasing\\order};
\end{tikzpicture}%
}

\newcommand{\RootedFigure}{%
\begin{tikzpicture}[>=Stealth,font=\small,x=1cm,y=.8cm]
\node[anchor=west] at (0,3.7){Fixed spacing: $D\le1$, $\sigma_a=\sigma_\infty$ for $a\ge1$};
\foreach \x in {1,3,5}{\fill (\x,2.8) circle (2pt);\draw[->] (\x,2.8)--(\x,3.4);}
\draw[->] (0,2.8)--(1,2.8);\draw[->] (1,2.8)--node[below]{$\sigma_\infty$}(3,2.8);
\draw[->] (3,2.8)--node[below]{$\sigma_\infty$}(5,2.8);\draw[->] (5,2.8)--(6,2.8);
\node[anchor=west] at (0,1.65){$K=10,m=4$: refinement $3\ (\alpha_3=2)\longrightarrow1\ (\alpha_1=1)$};
\node[anchor=west] at (0,.95){Spacing by age: $\sigma_0=2,\quad\sigma_1=4,\quad\sigma_2=\sigma_3=\cdots=5$};
\foreach \x in {1,3,5}{\fill (\x,0) circle (2pt);\draw[->] (\x,0)--(\x,.6);}
\draw[->] (0,0)--(1,0);\draw[->] (1,0)--node[below]{5}(3,0);\draw[->] (3,0)--node[below]{4}(5,0);\draw[->] (5,0)--(6,0);
\node[anchor=west] at (6.2,0){four copies};
\node[anchor=west] at (0,-1){Labels count cover edges; drawn segment lengths are schematic.};
\end{tikzpicture}}

\newcommand{\ContractionFigure}{%
\begin{tikzpicture}[>=Stealth,font=\small,x=.9cm,y=.65cm]
\draw[dotted] (0,0)--(3.2,0);
\draw[->] (0,0)--node[above left]{$U$}(1,1);
\draw[->] (1,1)--node[left=3pt]{$D$}(2,-1);
\draw[->] (2,-1)--node[below right]{$U$}(3,0);
\node[left=4pt] at (0,0){$3n+2$};
\node at (1.5,-1.9){weight $3n+2$};
\draw[->,double] (3.8,0)--(4.8,0);
\node (a) at (5.5,0){$n$};\node (b) at (8,0){$n$};
\draw[->] (a)--node[above]{$3n+2$}(b);
\node at (6.75,-1.15){one contribution to $Q_{n,n}$};
\node[draw,inner sep=5pt] (Q) at (1,-3.6){four-diagonal $Q$};
\node[draw,inner sep=5pt] (J) at (7.5,-3.6){tridiagonal $J$};
\draw[->] (Q)--node[above=3pt]{$QR=RJ$}(J);
\node at (4.25,-4.7){algebraic intertwining};
\end{tikzpicture}}

\title[Regular Pattern Packings and Moment Sequences]{%
  Permutations, Pattern Covers, Posets and Positivity, 
  Part II: Regular Pattern Packings and Moment Sequences}

\author{Natasha Blitvi\'c}
\address{School of Mathematical Sciences,
  Queen Mary University of London}
\email{n.blitvic@qmul.ac.uk}

\author{Mohammed Slim Kammoun}
\address{Universit\'e de Poitiers, CNRS, LMA}
\email{medslimkammoun@gmail.com}

\author{Einar Steingr\'imsson}
\address{Department of Mathematics and Statistics,
  University of Strathclyde}
\email{einar@alum.mit.edu}

\thanks{N. Blitvi\'c was partially supported by
  Leverhulme Trust Research Project Grant RPG-2020-103
  and the Engineering and Physical Sciences Research
  Council (EPSRC) grant EP/V048902.}

\thanks{M. S. Kammoun was partially supported by
  Leverhulme Trust Research Project Grant RPG-2020-103,
  CNRS-INSMI PEPS JCJC projects LAMPE and CAPRICES,
  LabEx CIMI (ANR-11-LABX-0040), and ANR project
  LOUCCOUM (ANR-24-CE40-7809).}

\date{\today}

\hypersetup{
  pdftitle={Permutations, Pattern Covers, Posets and Positivity
    Part II:  Moment Sequences},
  pdfauthor={Natasha Blitvic, Mohammed Slim Kammoun,
    Einar Steingrimsson}
}

\begin{document}

\begin{abstract}
We prove that broad families of regular pattern packings are enumerated
by moments of probability measures on $[0,\infty)$. A packing prescribes
occurrences of a consecutive permutation pattern at equally spaced
positions, with successive occurrences overlapping and together covering
the permutation. This regularity underlies the structure theory for the
associated pattern-packing posets that we developed in Part~I \cite{BKS1}. 

We obtain moment representations for all two-level shifted-grid
families of Part~I and all overlap-one patterns whose first and last
letters differ by one or two. For  families whose posets are trees with a unique minimum or maximum
and overlap at most half the pattern length, positivity
instead holds after finitely many initial terms of the enumerating sequence have been discarded. Our proofs develop further enumeration
formulas and combine probabilistic constructions with total
nonnegativity and weighted paths, yielding explicit laws or continued
fractions in several cases.

Although a previously conjectured universal moment property fails,
our results establish positivity for all packing sequences associated
with patterns of length at most four and all but one of the distinct
sequences arising from patterns of length five. Extensive numerical
evidence suggests that the remaining sequence is also a moment
sequence. These results lead us to formulate refined conjectures
for all overlap-one patterns and for the shifted-grid families of
Part~I with arbitrarily many levels.
\end{abstract}
\maketitle
\thispagestyle{empty}

\section{Introduction}

The theory of permutation patterns is rich in questions whose elementary
formulations give little indication of their difficulty. An occurrence
of a pattern $\pi$ in a permutation $\sigma$ is a subsequence whose
letters have the same relative order as those of $\pi$. Thus $4825$
in $64178235$ is an occurrence of $2413$, and a permutation avoids
$\pi$ if it contains no such subsequence. Even the growth rate of
permutations avoiding the single pattern $1324$, let alone their enumeration, remains a longstanding
open problem (see~\cite{brignall2026} and~\cite[Section~4]{BostanEtAl2020}). In addition to direct enumeration problems and growth rate analysis, questions about which patterns have equinumerous  avoidance have likewise led to classification problems that have occupied researchers for decades. The subject has meanwhile expanded in many directions, from patterns with positional restrictions to the distribution of occurrences in random permutations, with connections to sorting algorithms, partially ordered sets, topology, and algebraic combinatorics.

The theory of moments offers a perspective on such questions that brings
combinatorial sequences into contact with probability and analysis.
The factorials and Bell numbers are the moments of the exponential and
Poisson laws of mean one, respectively, while the Catalan numbers are
the even moments of the standard semicircle law. These familiar
correspondences belong to a broader collection of connections between
combinatorial structure and probability measures. The continued-fraction
framework of~\cite{BlitvicSteingrimsson2021}, for example, gives a
systematic interpretation of a large family of classical and
noncommutative probability laws in terms of permutation statistics.
Sokal and Zeng~\cite{SokalZeng2022} develop another broad
continued-fraction framework, refining enumerations of permutations,
set partitions, and perfect matchings by many simultaneous statistics.
For pattern-avoidance sequences, Bostan, Elvey Price, Guttmann and
Maillard~\cite{BostanEtAl2020} demonstrate how representing measures
can provide information about growth, including conditional bounds for
$1324$-avoiders under a conjectured moment property. The breadth of
these connections suggests that the moment property may reflect
features of combinatorial structure extending beyond any particular
enumeration. The problem is to identify those features and to determine
when they suffice to produce a representing measure.

We pursue this question for permutations with prescribed, regularly
spaced occurrences of a consecutive pattern. An occurrence is
\emph{consecutive} when its letters occupy adjacent positions.
Following Part~I~\cite{BKS1}, fix a pattern $\pi$ of length $K$ and an
overlap $1\le m<K$. A permutation \emph{$m$-packs $j$ occurrences of
$\pi$} if it has length $K+(j-1)(K-m)$ and contains occurrences
starting at
\[
1,\quad 1+(K-m),\quad 1+2(K-m),\quad\ldots,\quad
1+(j-1)(K-m).
\]
Successive prescribed occurrences therefore overlap in exactly $m$
positions and together cover the permutation. Additional occurrences
are permitted. We denote the number of such permutations by
$N_\pi^m(j)$ and consider admissible overlaps, for which the first $m$
letters of $\pi$ are order isomorphic to its last $m$ letters.
The definitions and notation are recalled in
Section~\ref{sec:preliminaries}.

The distinguishing feature of this construction is its \emph{regularity}.
The study of overlapping consecutive occurrences has a substantial
history, notably through cluster methods~\cite{ElizaldeNoy2012}, and
the individual numbers $N_\pi^m(j)$ are instances of the refined
cluster numbers of Dwyer and Elizalde~\cite[Section~2.1]{DwyerElizalde2018},
which prescribe overlapping occurrences without requiring equal spacing.
Our approach is to hold both the pattern and the overlap fixed as the
number of occurrences varies. Each additional occurrence then repeats
the same rule for identifying positions, making it possible to study
how fixed features of $\pi$ and $m$ govern an entire family of order
structures.

Each occurrence imposes a chain of inequalities on the values assigned
to its positions. Together these chains define the
\emph{pattern-packing poset} $\mathcal P_\pi^m(j)$, whose linear
extensions are in bijection with the packing permutations
\cite[Definition~5.2 and Theorem~5.3]{BKS1}. In Part~I, we developed
a structure theory for these posets. We described their branches,
which are chains attached to the rest of the poset by a single edge,
characterized the tree and chain cases, and constructed families of
shifted grid posets~\cite[Section~5]{BKS1}. These structural descriptions
are uniform in the number of occurrences. They allow the moment
problem to be approached through properties shared by an entire
family of posets, and form the starting point for the present paper.

More precisely, for fixed $\pi$ and $m$, write
\[
A_r=N_\pi^m(r+1),\qquad r\ge0.
\]
We seek a positive Borel measure $\mu$ on $\mathbb R$ such that
\begin{equation}\label{eq:intro-moments}
A_r=\int_{\mathbb R}x^r\,\mu(dx),\qquad r\ge0.
\end{equation}
Since $A_0=1$, any representing measure is a probability measure.
In the terminology of moment theory, $(A_r)_{r\ge0}$ is then a
\emph{Hamburger moment sequence}, and a \emph{Stieltjes moment sequence}
if a representing measure can be supported on $[0,\infty)$.
Hamburger's theorem characterizes~\eqref{eq:intro-moments} by positive
semidefiniteness of every Hankel matrix $(A_{i+j})_{0\le i,j\le n}$.
The Stieltjes criterion imposes the same condition on the shifted
matrices $(A_{i+j+1})_{0\le i,j\le n}$ as well.
The existence of a representing measure thus entails infinitely many
inequalities linking counts for different numbers of occurrences.
Establishing these inequalities is a separate problem from evaluating
the individual counts, and calls for further use of the underlying
combinatorial structure.

In this paper, we use the structure theory of Part~I to establish
moment positivity for broad families of regular pattern packings,
deriving further enumeration formulas in the course of the proofs.
The representing measures are supported on $[0,\infty)$ and are
explicit in several cases. For some tree families, positivity
requires the removal of a finite initial segment from the enumerating sequence, and we bound
the number of terms that must be discarded.

Our results establish positivity from the first term for every pattern
of length at most four and every admissible overlap
(Corollary~\ref{cor:small}), and for thirteen of the fourteen distinct
sequences arising from patterns of length at most five in the
motivating table of
\cite[Table~1]{BlitvicKammounSteingrimsson2023}.
Table~\ref{tab:moment-status} summarizes these results.
Extensive numerical evidence supports positivity of the remaining
sequence, shared by $13524$ and $24513$ at overlap two, but its
moment property remains unproved. The conjecture of
\cite[Conjecture~1]{BlitvicKammounSteingrimsson2023} that every
packing sequence admits a representing measure on $\mathbb R$
nevertheless fails.
We formulate refined conjectures asserting positivity for all
overlap-one patterns and for the shifted-grid families of Part~I
with any number of levels
(Conjectures~\ref{conj:overlap-one} and~\ref{conj:canonical-grids}).

The paper is organized as follows.
Section~\ref{sec:preliminaries} recalls the definitions and the
required results from Part~I. Sections~\ref{sec:enumeration}
and~\ref{sec:moments} derive enumeration formulas and moment
representations for overlap-one patterns with
$|\pi_K-\pi_1|=1$, beginning with the example $2413$.
Section~\ref{sec:products} brings together cases in which the
enumeration formulas factor directly into moments of independent
random variables. Section~\ref{sec:rooted} treats tree posets with a unique minimum
or maximum, assuming $m\le K/2$, and constructs moment
representations after a bounded number of initial terms have
been removed. Section~\ref{sec:grids} establishes positivity
for the two-level shifted-grid families of Part~I.
Section~\ref{sec:nonrooted} returns to overlap one and constructs
representing measures when $|\pi_K-\pi_1|=2$ and neither endpoint
is $1$ or $K$. Together with the earlier results, this establishes
positivity for all overlap-one patterns with
$|\pi_K-\pi_1|\le2$. Section~\ref{sec:explicit} constructs a compactly supported law
for $135246$ at overlap three, determines the Jacobi continued
fraction for $14325$ at overlap two, and obtains branched continued
fractions for a larger family of patterns.
Section~\ref{sec:scope} draws together the short-pattern results,
a counterexample to universal positivity, and the refined
conjectures.

\begin{table}[htbp]
\centering
\small
\setlength{\tabcolsep}{7pt}
\begin{tabular}{@{}lccc@{}}
\toprule
Pattern $\pi$ & $m=1$ & $m=2$ & $m=3$ \\
\midrule
$132$ & Prop.~\ref{prop:extreme-endpoint} & $-$ & $-$ \\
\addlinespace[4pt]
$1243$ & Prop.~\ref{prop:extreme-endpoint} & $-$ & $-$ \\
$1324$ & Prop.~\ref{prop:extreme-endpoint} & Cor.~\ref{cor:fuss} & $-$ \\
$1342$ & Prop.~\ref{prop:extreme-endpoint} & $-$ & $-$ \\
$1423$ & Prop.~\ref{prop:extreme-endpoint} & Cor.~\ref{cor:small} & $-$ \\
$2143$ & Thm.~\ref{thm:moments} & Cor.~\ref{cor:small} & $-$ \\
$2413$ & Thm.~\ref{thm:moments} & Cor.~\ref{cor:fuss} & $-$ \\
\addlinespace[4pt]
$12354$ & Prop.~\ref{prop:extreme-endpoint} & $-$ & $-$ \\
$12435$ & Prop.~\ref{prop:extreme-endpoint} & Thm.~\ref{thm:rooted-tails} & $-$ \\
$12453$ & Prop.~\ref{prop:extreme-endpoint} & $-$ & $-$ \\
$12534$ & Prop.~\ref{prop:extreme-endpoint} & Thm.~\ref{thm:rooted-tails} & $-$ \\
$13254$ & Prop.~\ref{prop:extreme-endpoint} & $-$ & Prop.~\ref{prop:short-overlap-three} \\
$13425$ & Prop.~\ref{prop:extreme-endpoint} & Cor.~\ref{cor:fuss} & $-$ \\
$13452$ & Prop.~\ref{prop:extreme-endpoint} & $-$ & $-$ \\
$13524$ & Prop.~\ref{prop:extreme-endpoint} & $\gamma_n$ & $-$ \\
$14253$ & Prop.~\ref{prop:extreme-endpoint} & $-$ & Prop.~\ref{prop:short-overlap-three} \\
$14325$ & Prop.~\ref{prop:extreme-endpoint} & Thm.~\ref{thm:grid} & $-$ \\
$14523$ & Prop.~\ref{prop:extreme-endpoint} & Thm.~\ref{thm:rooted-tails} & $-$ \\
$15234$ & Prop.~\ref{prop:extreme-endpoint} & Thm.~\ref{thm:rooted-tails} & $-$ \\
$15243$ & Prop.~\ref{prop:extreme-endpoint} & $-$ & Prop.~\ref{prop:short-overlap-three} \\
$15324$ & Prop.~\ref{prop:extreme-endpoint} & Thm.~\ref{thm:prefix} & $-$ \\
$15423$ & Prop.~\ref{prop:extreme-endpoint} & Thm.~\ref{thm:rooted-tails} & $-$ \\
$21354$ & Cor.~\ref{cor:21354} & Thm.~\ref{thm:rooted-tails} & $-$ \\
$21453$ & Thm.~\ref{thm:moments} & Thm.~\ref{thm:rooted-tails} & $-$ \\
$21534$ & Thm.~\ref{thm:gap-two} & $-$ & $-$ \\
$21543$ & Thm.~\ref{thm:moments} & Thm.~\ref{thm:rooted-tails} & $-$ \\
$23514$ & Thm.~\ref{thm:gap-two} & Thm.~\ref{thm:rooted-tails} & $-$ \\
$24153$ & Thm.~\ref{thm:moments} & $-$ & $-$ \\
$24513$ & Thm.~\ref{thm:moments} & $\gamma_n$ & $-$ \\
$25314$ & Thm.~\ref{thm:gap-two} & Thm.~\ref{thm:grid} & $-$ \\
$25413$ & Thm.~\ref{thm:moments} & Cor.~\ref{cor:fuss} & $-$ \\
\bottomrule
\end{tabular}
\caption{Positivity results for the patterns of \cite[Table~1]{BKS1}. One sequence, denoted $(\gamma_n)$ and appearing in two entries, remains open. Dashes denote inadmissible overlaps. }
\label{tab:moment-status}
\end{table}


\section{Preliminaries: pattern packings}\label{sec:preliminaries}

We recall the definition of pattern packing from
Part~I~\cite[Definition~2.6 and Remark~2.7]{BKS1}. As is habitual in combinatorics, we let $[n]=\{1,\ldots,n\}$ and denote by $S_n$ the set of all permutations of $[n]$, that is, the bijections $[n]\to[n]$. We interpret elements of $S_n$ as words of length $n$ in the alphabet $[n]$ formed of distinct letters, and write $\sigma=\sigma_1\cdots\sigma_n$. Two words of distinct real numbers are order isomorphic,
written $w\sim w'$, if their letters have the same relative
order. For example, $4936\sim2413$.

A consecutive occurrence of $\pi\in S_K$ in $\sigma\in S_n$
starting at position $i$ is a block
$\sigma_i\cdots\sigma_{i+K-1}$ order isomorphic to $\pi$.
We refer to $\pi$ as a $K$-pattern.

\begin{definition}[Pattern packing]\label{def:packing}
Fix $K\ge2$, $j\ge1$, a pattern $\pi\in S_K$, and $m<K$, with $m$ positive. Set
\[
d=K-m,\qquad n_j=K+(j-1)d.
\]
A permutation $\sigma\in S_{n_j}$ \emph{$m$-packs $j$ occurrences of $\pi$} if
\begin{equation}\label{eq:packing}
\sigma_{1+td}\cdots\sigma_{K+td}\sim\pi
\qquad (0\le t<j).
\end{equation}
Write $\mathcal S_\pi^m(j)$ for this set and $N_\pi^m(j)$ for its cardinality.
\end{definition}

Definition~\ref{def:packing} does not preclude additional occurrences outside the prescribed ones. For example, $192837465$ 1-packs two occurrences of $15243$, starting at positions 1 and 5, and also contains an occurrence starting at position 3. The parameter $m$ is the number of shared positions in consecutive prescribed occurrences, and consecutive prescribed occurrences begin $d=K-m$ positions apart.

An overlap $m$ is \emph{admissible} for $\pi$ if the $m$-prefix and $m$-suffix of $\pi$ are order isomorphic, that is:
\[
\pi_1\cdots\pi_m\sim\pi_{K-m+1}\cdots\pi_K.
\]
We assume admissibility whenever more than one occurrence is considered. In particular, overlap one is admissible for every pattern.
Following Part~I~\cite{BKS1}, the \emph{interfix} of a pattern $\pi$, given overlap $m$, is the word between the $m$-prefix and $m$-suffix when $m<K/2$; it is empty when $m\ge K/2$.

\subsection{The associated poset}

 The act of packing a given number of occurrences of a consecutive permutation pattern forces relations between the letters of the permutation. For example, any permutation $\sigma$ 1-packing two occurrences of $2413$ must simultaneously satisfy
\begin{equation}\label{eq:2413-two}
\sigma_3<\sigma_1<\sigma_4<\sigma_2,
\qquad
\sigma_6<\sigma_4<\sigma_7<\sigma_5.
\end{equation}
This implies, for example, that $\sigma_2>\sigma_6$, because of the interaction between the two occurrences at their shared position.

In general, given a $K$-pattern $\pi$, let $\tau=\pi^{-1}$, so that $\tau_i$ is the position of the entry $i$ in $\pi$.
The \emph{pattern-packing poset} $\ppmj$ has ground
set $[n_j]$ and its order is the transitive closure of the strict orders generated by the chains
\begin{equation}\label{eq:poset-chains}
td+\tau_1< td+\tau_2<\cdots< td+\tau_K
\qquad(0\le t<j).
\end{equation}
(See \cite[Definition~5.2]{BKS1}.) These chains are the lines of the \emph{stack} in Part~I~\cite{BKS1}. We also refer to each as the chain of its occurrence.
The elements of $\ppmj$ represent positions in a $\pi$-packing
permutation, and the order specifies when a position $x$ must have a
smaller value than a position $y$.

For an admissible overlap, these relations are consistent.
Indeed, each new occurrence meets the existing poset precisely
in the $m$ positions it shares with the preceding occurrence.
These positions form a chain in both copies, with the same
ordering by admissibility. Any path through newly added vertices
between two old vertices therefore imposes a comparison that
already holds in the existing poset. Consequently, the attachment
cannot create a directed cycle. Induction on the number of
occurrences shows that the construction defines a partial order.

A \emph{linear extension} of a poset is a total ordering of its elements in which $u$ precedes $v$ whenever $u< v$. If $v_1,\ldots,v_{n_j}$ is a linear extension of $\ppmj$, assign the value $i$ to position $v_i$. The resulting permutation satisfies~\eqref{eq:packing}. Conversely, listing the positions of a packing permutation in increasing order of their values gives that linear extension. Thus, as in Part~I~\cite[Theorem~5.3]{BKS1},
\begin{equation}\label{eq:extensions}
N_\pi^m(j)=|\mathcal L(\ppmj)|.
\end{equation}
As in \cite[Corollary~5.21]{BKS1}, we call a packing poset a \emph{tree poset}, or just a tree, when its undirected Hasse graph is a tree. In particular, for overlap one the poset is obtained by attaching successive occurrence chains at a single vertex and is a tree. 
We also refer to the \emph{core} and \emph{branches} of Part~I~\cite{BKS1}. The core is the interval between the designated bottom and top elements. A branch is a nonempty maximal connected part outside the core, which by \cite[Theorem~5.16]{BKS1} is a chain attached to the core by a single edge. Its length counts its elements, excluding the attachment vertex in the core. When the poset is a tree, the core is itself a chain.

By the equivalences in~\cite[Section~3]{BKS1}, when $m\le K/2$, the enumeration for a pattern $\pi$ is invariant under reversal and complementation of $\pi$, as well as by sorting the prefix and suffix, exchanging them, permuting the interfix arbitrarily, and by any composition of these operations. We therefore typically, when $m\le K/2$, consider only patterns with increasing prefix, suffix and interfix, and for which the first letter of the prefix is smaller than the first letter of the suffix.

\section{Tree enumeration when $|\pi_K-\pi_1|=1$}\label{sec:enumeration}

We first consider overlap $m=1$ and $K$-patterns whose first and
last letters have consecutive values, $|\pi_K-\pi_1|=1$.
Since reversing both the pattern and the packed permutation
preserves the counts, we may assume that $\pi_K=\pi_1+1$.

\subsection{Order statistics and the density method}\label{sec:labels}

We start by reformulating the problem of counting linear extensions in terms of order statistics of random variables.
Let $P$ be a poset on $[n]$, and consider independent,
identically distributed real random variables
$X_1,\ldots,X_n$ with a common distribution without atoms.
Their values are almost surely distinct. By exchangeability,
each of the $n!$ possible orderings of these variables is
equally likely. Equivalently, their order statistics
$X_{(1)}<\cdots<X_{(n)}$ can be assigned to the vertices
of $P$ in $n!$ equally likely ways.

The event that the random variables respect the order of
the poset is
\[
E_P=\{X_u<X_v\text{ whenever }u<_P v\}.
\]
Since each linear extension corresponds to one of these
orderings, we have
\begin{equation}\label{eq:order-probability}
|\mathcal L(P)|=n!\,\Prob(E_P).
\end{equation}
For a packing poset $P=\ppmj$, the event $E_P$
is, almost surely, the event that the word
$X_1\cdots X_{n_j}$ is order isomorphic to a permutation
in $\mathcal S_\pi^m(j)$.

For uniform labels on $[0,1]$, this probability is the
volume of the region of the unit cube defined by the order
inequalities, namely the order polytope of $P$ (up to its
boundary, which has measure zero).
This interpretation underlies the \emph{density method}
of Banderier and Wallner~\cite[Sec.~3]{BanderierWallner2021},
who assign independent uniform labels to cells and compute
the volume of the region satisfying the ordering constraints
by conditioning on selected coordinates and integrating
over the others.
While the order-statistic interpretation applies to any
common atomless distribution, and the corresponding
integral representations may take different forms,
the ordering probabilities are necessarily the same.
Indeed, if $F$ is the common distribution function,
then $U_i=F(X_i)$ are independent uniform random variables
whose relative ordering agrees almost surely with that
of the original variables.

\subsection{The example 2413}\label{sec:2413}

We next show how to use order statistics of independent and identically distributed random variables to compute the number of permutations 1-packing $j$ occurrences of $2413$. This example exhibits both conditional independence and the symmetry that reduces a multiple integral to a power of a single integral.

\begin{figure}[htbp]
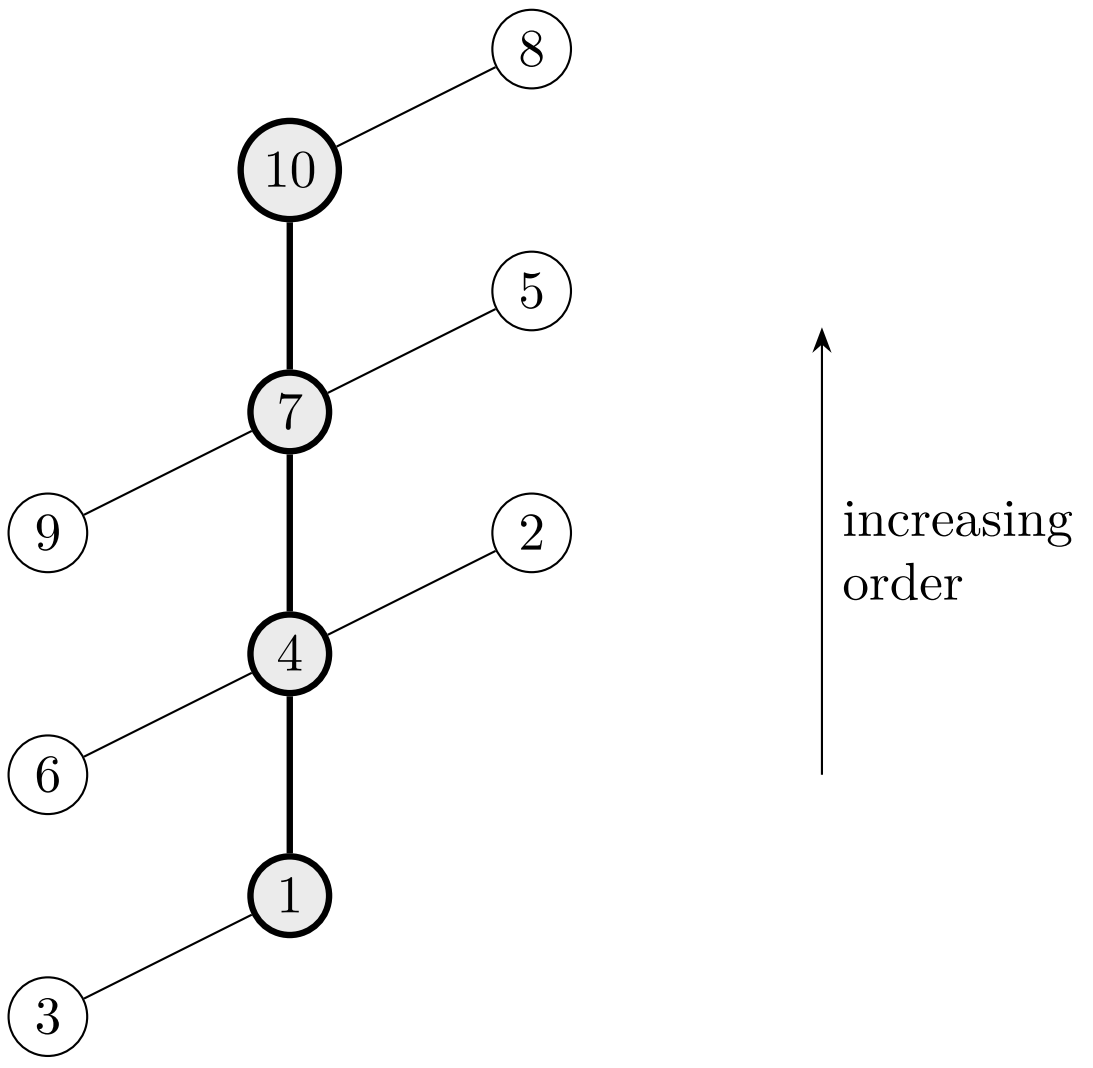

\centering\PackingTreeFigure
\caption{The Hasse diagram for three overlap-one occurrences of $2413$. The vertices are positions, with order increasing upward. Each occurrence contributes the chain $3+3t<_P1+3t<_P4+3t<_P2+3t$. The highlighted trunk consists of positions $1,4,7,10$. The remaining labels are independent after the trunk labels have been fixed.}\label{fig:2413}
\end{figure}

Let $X_1,\ldots,X_{1+3j}$ be independent random variables with common continuous cumulative distribution function $F$. In occurrence $t$, where $0\le t<j$, the required inequalities are
\begin{equation}\label{eq:2413-labels}
X_{3+3t}<X_{1+3t}<X_{4+3t}<X_{2+3t}.
\end{equation}
The labels $X_1,X_4,\ldots,X_{1+3j}$ lie on the common increasing path in the tree. We call these the \emph{trunk variables}. The trunk is a subchain of the core, and the remaining variables belong to branches or to the initial and final portions of the core. Figure~\ref{fig:2413} displays the construction for three occurrences.

By conditioning on the trunk variables, two observations follow. First, with $X_{1+3t}=y_t$ and $y_0<\cdots<y_j$, the conditions on the remaining variables are
\[
X_{3+3t}<y_t,\qquad X_{2+3t}>y_{t+1}\qquad(0\le t<j).
\]
These events involve distinct independent labels. Their conditional probability is therefore
\begin{equation}\label{eq:2413-conditional}
\prod_{t=0}^{j-1}F(y_t)(1-F(y_{t+1}))
=F(y_0)(1-F(y_j))\prod_{t=1}^{j-1}F(y_t)(1-F(y_t)).
\end{equation}

Second, with the values $y_0=y$ and $y_j=z$ fixed, the product involving the internal trunk variables is symmetric in those variables. The product measure $dF(y_1)\cdots dF(y_{j-1})$ gives zero mass to coincidences. Its integral over $y<y_1<\cdots<y_{j-1}<z$ is therefore $1/(j-1)!$ times its integral over the whole cube $(y,z)^{j-1}$:
\begin{align}\label{eq:2413-simplex}
&\int_{y<y_1<\cdots<y_{j-1}<z}
\prod_{t=1}^{j-1}F(y_t)(1-F(y_t))\,dF(y_t)\notag\\
&\hspace{30mm}=\frac{1}{(j-1)!}
\left(\int_{(y,z)}F(x)(1-F(x))\,dF(x)\right)^{j-1}.
\end{align}
Finally, taking expectations over the bottom and top trunk labels and using~\eqref{eq:order-probability} gives
\begin{align}\label{eq:2413-F}
N_{2413}^1(j)
={}&\frac{(1+3j)!}{(j-1)!}
\int_{y<z}F(y)(1-F(z))\notag\\
&\hspace{19mm}\cdot
\left(\int_{(y,z)}F(x)(1-F(x))\,dF(x)\right)^{j-1}
\,dF(y)\,dF(z).
\end{align}
For $j=1$, the internal product and integral over zero coordinates in~\eqref{eq:2413-simplex} are both interpreted as one.

In particular, taking the labels to be uniform on $[0,1]$ gives
\begin{equation}\label{eq:2413-uniform}
N_{2413}^1(j)=\frac{(1+3j)!}{(j-1)!}
\int_0^1\int_y^1 y(1-z)
\left(\frac{z^2-y^2}{2}-\frac{z^3-y^3}{3}\right)^{j-1}dz\,dy.
\end{equation}
Here the expression in parentheses is $\int_y^z x(1-x)\,dx$. The example has reduced the enumeration to the moments of a function of two endpoint variables, together with a factorial factor.

\subsection{The general integral}\label{sec:general-integral}

Let $\pi$ be a $K$-pattern with $\pi_K=\pi_1+1$, and set
\begin{equation}\label{eq:parameters}
a=\pi_1-1,\qquad b=K-\pi_K,\qquad d=K-1.
\end{equation}
Since $K=a+b+2$, the remaining positions in each
occurrence form two ordered chains: one consisting of
$a$ positions below the first endpoint, and the other
of $b$ positions above the last endpoint.
Either chain may be empty.

\begin{proposition}[Adjacent-endpoint enumeration]\label{prop:integral}
Let $\pi\in S_K$, $K\ge2$, satisfy $\pi_K=\pi_1+1$, and use~\eqref{eq:parameters}. For every $j\ge1$,
\begin{equation}\label{eq:adjacent-count}
N_\pi^1(j)=\frac{(1+jd)!}{(j-1)!(a!b!)^j}
\int_{0<y<z<1}y^a(1-z)^b
\left(\int_y^z x^a(1-x)^b\,dx\right)^{j-1}dy\,dz.
\end{equation}
\end{proposition}
\begin{remark}
Equivalently, for any continuous cumulative distribution
function $F$, the integral in~\eqref{eq:adjacent-count}
admits the representation
\begin{equation}\label{eq:adjacent-F}
\int_{y<z}F(y)^a(1-F(z))^b
\left(\int_{(y,z)}F(x)^a(1-F(x))^b\,dF(x)\right)^{j-1}
\,dF(y)\,dF(z).
\end{equation}
\end{remark}
\begin{proof}
The proof follows the same approach as the calculation
for $2413$, with ordered chains now playing the role of
the individual variables outside the trunk. As before, we integrate
out these variables while retaining the trunk labels,
and then exploit the symmetry of the internal trunk
coordinates.

Let $X_1,\ldots,X_{1+jd}$ be independent random variables
with common continuous cumulative distribution function $F$, and set
$T_t=X_{1+td}$ for $0\le t\le j$.
The ordering constraints along the trunk require
$T_0<T_1<\cdots<T_j$.

Condition on the trunk labels taking values
$y_0<\cdots<y_j$.
For each occurrence $t$, the $a$ variables on the lower chain
must all lie below $y_t$, in a prescribed relative order.
The probability of this event is $F(y_t)^a/a!$:
the probability that all $a$ labels lie below $y_t$
is $F(y_t)^a$, and by exchangeability, each of their
$a!$ possible relative orders is equally likely.
Similarly, the probability that the $b$ variables on the upper
chain lie above $y_{t+1}$ in the required order is
$(1-F(y_{t+1}))^b/b!$.

Since these lower and upper chains occupy disjoint positions in distinct occurrences,
these events are independent conditional on the trunk labels.
Their joint probability is therefore
\begin{align}\label{eq:general-conditional}
\prod_{t=0}^{j-1}
\frac{F(y_t)^a(1-F(y_{t+1}))^b}{a!b!}
&=(a!b!)^{-j}F(y_0)^a(1-F(y_j))^b
\prod_{t=1}^{j-1}F(y_t)^a(1-F(y_t))^b.
\end{align}
The same expression holds when $a=0$ or $b=0$,
with an empty chain contributing a factor of one.

Now fix the endpoint labels $y_0=y$ and $y_j=z$.
The product involving the remaining trunk labels
$y_1,\ldots,y_{j-1}$ is symmetric in these variables.
As in~\eqref{eq:2413-simplex}, integration over the
ordered simplex $y<y_1<\cdots<y_{j-1}<z$ therefore gives
\[
\frac{1}{(j-1)!}
\left(
\int_{(y,z)}F(x)^a(1-F(x))^b\,dF(x)
\right)^{j-1}.
\]
Integrating over the endpoint labels $T_0$ and $T_j$
then yields the required order probability.
Multiplying by $(1+jd)!$ gives the counting formula
with the prefactor in~\eqref{eq:adjacent-count}
and the integral in~\eqref{eq:adjacent-F}.
Taking $F$ to be the uniform distribution gives
\eqref{eq:adjacent-count}.
(Since the integrands are bounded and measurable
and the measures are finite, the interchanges
of integration are justified by Fubini's theorem.)
\end{proof}


\begin{remark}
The integral~\eqref{eq:adjacent-count}
reduces to the familiar, elementary cases.
For $j=1$, we have
\begin{equation}\label{eq:triangle-mass}
\int_0^1\int_0^z y^a(1-z)^b\,dy\,dz
=\frac{1}{a+1}\int_0^1 z^{a+1}(1-z)^b\,dz
=\frac{a!b!}{(a+b+2)!}
=\frac{a!b!}{K!},
\end{equation}
giving $N_\pi^1(1)=1$.

When $K=2$, we have $a=b=0$ and $d=1$, so
\[
\int_{0<y<z<1}(z-y)^{j-1}\,dy\,dz
=\frac{1}{j(j+1)}.
\]
Multiplication by $(j+1)!/(j-1)!$
again gives $N_\pi^1(j)=1$.
\end{remark}

\section{A moment representation when $|\pi_K-\pi_1|=1$}\label{sec:moments}

The integral in~\eqref{eq:adjacent-count} has an evocative form: its dependence on $j$
appears through the $(j-1)$st power of a nonnegative
function, integrated against a weight independent of $j$. To obtain a moment representation for
the counting sequence, however, we must also account for
the factor
$\frac{(1+(r+1)d)!}{r!(a!b!)^{r+1}}.$ For this, we use the following elementary lemma.

\begin{lemma}\label{lem:density-transform}
Let $q$ be a bounded nonnegative nonincreasing probability density on $(0,M)$, where $0<M<\infty$. Choose a right-continuous version and extend it by zero on $[M,\infty)$. Let $\rho=-dq$ be its positive Lebesgue--Stieltjes measure on $(0,\infty)$. Then
\[
\sigma(dt)=t\,\rho(dt)
\]
is a probability measure supported on $(0,M]$, and
\begin{equation}\label{eq:density-transform}
\int t^r\,\sigma(dt)=(r+1)\int_0^M s^r q(s)\,ds
\qquad(r\ge0).
\end{equation}
\end{lemma}

\begin{proof}
The monotonicity of $q$ gives a finite positive measure $\rho$ characterized by
\[
\rho((u,v])=q(u)-q(v)\qquad(0<u<v).
\]
In particular, $q(s)=\rho((s,\infty))$ for $s>0$, and $\rho$ is supported on $(0,M]$. Tonelli's theorem gives
\begin{align*}
\int_0^M s^r q(s)\,ds
&=\int_{(0,M]}\left(\int_0^t s^r\,ds\right)\rho(dt)\\
&=\frac{1}{r+1}\int_{(0,M]} t^{r+1}\,\rho(dt).
\end{align*}
This proves~\eqref{eq:density-transform}.
Taking $r=0$ and using the fact that $q$ is a
probability density shows that $\sigma$ has
total mass one.

When $q$ is absolutely continuous on $(0,M)$,
the measure $\sigma$ has density $-tq'(t)$
on this interval. Write
$q(M-)=\lim_{t\uparrow M}q(t)$.
If $q(M-)>0$, extending $q$ by zero at $M$
introduces a discontinuity, which contributes
an atom of mass $Mq(M-)$ to $\sigma$ at $M$.
\end{proof}

\begin{theorem}\label{thm:moments}
For every $K\ge2$ and every pattern $\pi\in S_K$ with $|\pi_K-\pi_1|=1$, the sequence
\[
\bigl(N_\pi^1(r+1)\bigr)_{r\ge0}
\]
is the moment sequence of a probability measure on $[0,\infty)$.
\end{theorem}

\begin{proof}
By reversal, assume $\pi_K=\pi_1+1$, and use the parameters in~\eqref{eq:parameters}. The case $K=2$ has already given the constant sequence one and is represented by the point mass at 1. Henceforth $d=K-1\ge2$.

The proof proceeds in three steps. We first interpret
the integral in~\eqref{eq:adjacent-count} as the moments
of a compactly supported random variable. This leaves
the prefactor, which depends on $r$, to be accounted for.
We do so by using a decreasing-density transformation
to incorporate the factor $r+1$, and then representing
the remaining factorial ratio as the moments of a
product of independent gamma random variables.

\medskip
\noindent\emph{(1) The first moment integral.}
We analyze the sequence $(I_r)$  defined by
\[
I_r=\int_0^1\int_0^z y^a(1-z)^b
\left(\int_y^z x^a(1-x)^b\,dx\right)^rdy\,dz.
\]
Let
\begin{equation}\label{eq:f-definition}
h(t)=\int_0^t x^a(1-x)^b\,dx,\qquad
M=h(1)=\frac{a!b!}{(a+b+1)!},\qquad c=\frac{K!}{a!b!}.
\end{equation}
The function $h$ is strictly increasing on $[0,1]$, so $h:[0,1]\to[0,M]$ is a bijection. Write $g=h^{-1}$. By~\eqref{eq:triangle-mass},
\[
p(y,z)=c\,y^a(1-z)^b\ind_{\{0<y<z<1\}}
\]
is a probability density. Let $(Y,Z)$ have this density and define
\begin{equation}\label{eq:X-definition}
X=\int_Y^Z x^a(1-x)^b\,dx=h(Z)-h(Y).
\end{equation}
By definition of expectation,
\begin{equation}\label{eq:X-moments}
\E X^r=cI_r\qquad(r\ge0).
\end{equation}
Moreover, $0<X<M$ almost surely. Its distribution has support $[0,M]$, since the joint density is positive throughout the open triangle and $h(z)-h(y)$ takes every value in $(0,M)$ there.

To determine the density of $X$, extend $g$ by $g(t)=1$ for $t>M$ and $g(t)=0$ for $t<0$. Since $h$ is increasing, the condition $h(Z)-h(Y)\le s$ is equivalent to $Z\le g(h(Y)+s)$ for $0\le s\le M$. Hence the cumulative distribution function is
\begin{equation}\label{eq:X-cdf}
\Prob(X\le s)
=c\int_0^1 y^a
\left(\int_y^{\min\{1,g(h(y)+s)\}}(1-z)^b\,dz\right)dy.
\end{equation}
The inverse satisfies
\[
g'(u)=\frac{1}{g(u)^a(1-g(u))^b},\qquad 0<u<M.
\]
For a fixed $y\in(0,1)$, make the substitution $u=h(z)-h(y)$ in the inner integral of~\eqref{eq:X-cdf}. The factor $(1-z)^b$ cancels the same factor in $h'(z)=z^a(1-z)^b$. Therefore
\begin{align*}
\Prob(X\le s)
&=c\int_0^1\int_0^{\min\{s,M-h(y)\}}
\left(\frac{y}{g(h(y)+u)}\right)^a du\,dy\\
&=\int_0^s c\int_0^{g(M-u)}
\left(\frac{y}{g(h(y)+u)}\right)^a dy\,du.
\end{align*}
The substitution can first be made with $z$ in a closed subinterval of $(y,1)$ and then extended to the endpoint by monotone convergence. The interchange is justified by Tonelli. It follows that $X$ has density
\begin{equation}\label{eq:q}
q(s)=c\int_0^{g(M-s)}
\left(\frac{y}{g(h(y)+s)}\right)^a dy,
\qquad 0<s<M.
\end{equation}

\medskip
\noindent\emph{(2) Incorporating  the factor $r+1$ via monotonicity.}
As $s$ increases, the upper integration limit $g(M-s)$ decreases. On the common part of the integration intervals, the denominator $g(h(y)+s)$ increases, so the nonnegative integrand in~\eqref{eq:q} is nonincreasing. Thus $q$ is nonincreasing. Since $g(h(y)+s)\ge y$, the integrand is at most one. With the convention that its zeroth power is one when $a=0$, dominated convergence gives
\[
q(0+)=c,\qquad q(M-)=0.
\]
The same domination, with the changing integration interval represented by an indicator on $(0,1)$, proves continuity for $0<s<M$.

Lemma~\ref{lem:density-transform} now yields a probability measure
\begin{equation}\label{eq:sigma}
\sigma(dt)=-t\,dq(t)
\end{equation}
on $[0,M]$, with no atom contributed at $M$, whose moments are
\begin{equation}\label{eq:V-moments}
\E V^r=(r+1)\E X^r=(r+1)cI_r,
\end{equation}
where $V$ has law $\sigma$. This is the factor by which the original integral moments must be multiplied before treating the factorial ratio.

\medskip
\noindent\emph{(3) The remaining factorial factor.}
For $\alpha>0$, let $G_\alpha$ have the gamma density
\[
\frac{x^{\alpha-1}e^{-x}}{\Gamma(\alpha)}\ind_{\{x>0\}}.
\]
Its $r$th moment is $\Gamma(\alpha+r)/\Gamma(\alpha)=(\alpha)_r$, where $(\alpha)_r=\alpha(\alpha+1)\cdots(\alpha+r-1)$ and $(\alpha)_0=1$.

At $j=r+1$, Proposition~\ref{prop:integral} and~\eqref{eq:V-moments} give
\begin{align}\label{eq:factor-target}
N_\pi^1(r+1)
&=\frac{(dr+d+1)!}{r!(a!b!)^{r+1}}\,I_r=\frac{(dr+d+1)!}{(d+1)!(a!b!)^r(r+1)!}\,\E V^r,
\end{align}
because $c=(d+1)!/(a!b!)$. The quotient of factorials can be rewritten as
\begin{align}\label{eq:residue-product}
\frac{(dr+d+1)!}{(d+1)!}
&=\prod_{\ell=0}^{r-1}\prod_{t=2}^{d+1}(\ell d+d+t)
=d^{dr}\prod_{t=2}^{d+1}(1+t/d)_r.
\end{align}
For $r=0$ both sides equal one. For $r\ge1$, the first product contains each integer from $d+2$ to $dr+d+1$ exactly once. Since $d\ge2$, the index $t=d$ is present, and its factor is $(2)_r=(r+1)!$. Canceling it gives
\begin{equation}\label{eq:gamma-factor}
\frac{(dr+d+1)!}{(d+1)!(a!b!)^r(r+1)!}
=\left(\frac{d^d}{(a!b!)}\right)^r
\prod_{\substack{2\le t\le d+1\\t\ne d}}(1+t/d)_r.
\end{equation}

Choose independent gamma variables $G_{1+t/d}$ for the indices in~\eqref{eq:gamma-factor}, independently of $V$, and set
\begin{equation}\label{eq:final-law}
W=\frac{d^d}{(a!b!)}\,V
\prod_{\substack{2\le t\le d+1\\t\ne d}}G_{1+t/d}.
\end{equation}
All factors are nonnegative, their moments of every order are finite, and independence gives
\[
\E W^r
=\frac{(dr+d+1)!}{(d+1)!(a!b!)^r(r+1)!}\,\E V^r
=N_\pi^1(r+1)
\qquad(r\ge0).
\]
The equality at $r=0$ verifies the total mass. No step required $a$ or $b$ to be positive, so the cases of an empty lower or upper chain are included. Reversal gives the decreasing endpoint orientation and completes the proof.
\end{proof}

For example, for the pattern $2413$, we have $a=b=1$, $d=3$,
$c=24$, and $M=1/6$, with $h(t)=t^2/2-t^3/3$.
Thus~\eqref{eq:final-law} becomes
\[
W=27\,V\,G_{5/3}G_{7/3},
\]
where the factors are independent and $V$ has the
distribution defined by~\eqref{eq:q}--\eqref{eq:sigma}.

\begin{remark}
The representing measure for $2413$ is Stieltjes determinate.
More generally, the same holds whenever
$|\pi_K-\pi_1|=1$ and $d\le3$.

For $d\ge2$, the variable $V$ in~\eqref{eq:final-law}
is supported on $[0,M]$, where $M=a!b!/d!$.
Each of the $d-1$ gamma factors has $r$th moment
at most $(r+2)^r$. Consequently,
\[
N_\pi^1(r+1)
=\E W^r
\le \left(\frac{d^d}{d!}\right)^r
(r+2)^{(d-1)r}
\qquad(r\ge1).
\]
It follows that
\[
\bigl(N_\pi^1(r+1)\bigr)^{-1/(2r)}
\ge
\left(\frac{d!}{d^d}\right)^{1/2}
(r+2)^{-(d-1)/2}.
\]
For $d=2,3$, the sum of these lower bounds diverges.
By the Stieltjes Carleman criterion
\cite[Theorem~1.11]{ST43},
the representing measure is unique among measures
on $[0,\infty)$.

For $d=1$, the sequence is identically one and its
unique representing measure is $\delta_1$.
\end{remark}

\section{Products of moments and direct enumeration formulas}\label{sec:products}

In the case $|\pi_K-\pi_1|=1$, converting counts to moments required a change of density before the counting factorial could be represented probabilistically. In contrast, some packing formulas factor directly into moments of independent random variables. We bring these cases together in this section.

\subsection{Products and binomial factors}

We remind the reader of the following elementary fact:

\begin{lemma}\label{lem:product-moments}
If $(a_r)_{r\ge0}$ and $(b_r)_{r\ge0}$ are Hamburger moment sequences, then $(a_rb_r)_{r\ge0}$ is a Hamburger moment sequence. If both are Stieltjes, their product is Stieltjes. For any $C\in\R$, $(C^ra_r)_{r\ge0}$ is Hamburger, and it is Stieltjes when $C\ge0$ and $(a_r)$ is Stieltjes.
\end{lemma}
\begin{proof}
Suppose first that $a_0=b_0=1$. Let $X$ and $Y$ be independent random variables with moments $(a_r)$ and $(b_r)$, respectively. The product $XY$ has moments
\[
\E[(XY)^r]=\E[X^r]\E[Y^r]=a_rb_r.
\]
If the two laws are supported on $[0,\infty)$, so is their product law. 

For $a_0,b_0>0$, we apply the probability argument to $a_r/a_0$ and $b_r/b_0$, then multiply the
resulting measure by $a_0b_0$. If either $a_0$ or $b_0$ is zero, the corresponding positive measure is zero and the product sequence is represented by the zero measure. 

Finally, multiplication by $C^r$ can be similarly viewed as product with an independent random variable given by a point mass at $C$.
\end{proof}

Next, we show how to express certain products of binomial coefficients as gamma moments. Recall the gamma variables $G_\alpha$ and rising factorials introduced in Section~\ref{sec:moments}. For $h>0$,  integer $\ell\ge0$, and $r\ge0$, rearranging a finite product gives
\begin{align}\label{eq:binomial-gammas}
\prod_{i=0}^{r}\binom{ih+\ell}{\ell}
=\prod_{t=1}^{\ell}\prod_{i=1}^{r}\frac{ih+t}{t}
=\left(\frac{h^\ell}{\ell!}\right)^r
\prod_{t=1}^{\ell}(1+t/h)_r
=\E\left(\frac{h^\ell}{\ell!}\prod_{t=1}^{\ell}G_{1+t/h}\right)^r,
\end{align}
where the gamma random variables are independent.  We identify the empty product, when $\ell=0$, with the constant random variable one. When $h=0$, the original binomial product is again one and is treated separately.

The enumeration in the following proposition is the
overlap-one case of Part~I~\cite[Proposition~5.25]{BKS1}
and appeared earlier in
\cite[Proposition~1]{BlitvicKammounSteingrimsson2023}.
Here, we identify its representing measure.

\begin{proposition}\label{prop:extreme-endpoint}
Let $\pi\in S_K$, $K\ge2$, with $\pi_1=1$, and let $\ell=K-\pi_K$. Then
\[
N_\pi^1(j)=\prod_{i=0}^{j-1}\binom{i(K-1)+\ell}{\ell},\qquad j\ge1.
\]
Consequently $(N_\pi^1(r+1))_{r\ge0}$ is Stieltjes, with the law in~\eqref{eq:binomial-gammas} for $h=K-1$. Reversal and value complementation give the same conclusion whenever the first or last letter is $1$ or $K$.
\end{proposition}
\begin{proof}
Let $s=\pi_K$ and consider a packing of $j\ge2$ occurrences, of length $n=1+j(K-1)$. The shared letter at position $K$ is the least letter in the second occurrence. It is smaller than every later letter, by induction through the successive overlaps. In the first occurrence, exactly $s-1$ letters lie below it. Those letters and the shared letter are therefore the values $1,\ldots,s$, in their prescribed relative order.

The $\ell=K-s$ letters of the first occurrence above the shared letter form a chain. They may be assigned any $\ell$ values from $\{s+1,\ldots,n\}$. There are
\[
\binom{n-s}{\ell}=\binom{(j-1)(K-1)+\ell}{\ell}
\]
choices. The remaining letters, together with the shared letter, form a packing of $j-1$ occurrences after standardization. Conversely, any such packing and any choice of the $\ell$ values reconstruct the original word uniquely. 
The first occurrence meets the remaining packing only at the shared vertex, which lies below both these $\ell$ vertices and all subsequently added vertices. Thus no vertex in the former group is comparable with any vertex in the latter.
Hence
\[
N_\pi^1(j)=\binom{(j-1)(K-1)+\ell}{\ell}N_\pi^1(j-1).
\]
Starting with $N_\pi^1(1)=1$ proves the product formula. Equation~\eqref{eq:binomial-gammas} proves the moment assertion. Reversal preserves the prescribed windows, and replacing every value $x$ by $n+1-x$ reverses all comparisons. Both operations preserve the number of packings.
\end{proof}

\subsection{Separated and coincident occurrences}

The packing definition for a pattern $\pi$ also has a direct extension to $m\le0$: consecutive prescribed windows are separated by $-m$ unconstrained positions. At $m=0$ they are disjoint and consecutive. At the other boundary, $m=K$, all prescribed windows coincide and there is only one permutation 
covered by $\pi$, namely $\pi$ itself. These conventions include $K=1$.

\begin{proposition}\label{prop:disjoint}
For $K\ge1$, an integer $m\le0$, and $d=K-m$,
\begin{equation}\label{eq:separated}
N_\pi^m(j)=\frac{(K+(j-1)d)!}{(K!)^j},\qquad
N_\pi^m(r+1)=\E\left(\frac{G_{K+1}^{\,d}}{K!}\right)^r.
\end{equation}
For $m=K$, the counts are all one and the representing law is $\delta_1$.
\end{proposition}
\begin{proof}
Since we do not prohibit further occurrences beyond the prescribed ones, each of the $j$ disjoint windows imposes only its own relative order. The packing poset is a disjoint union of $j$ chains of length $K$ and $(-m)(j-1)$ isolated vertices. For independent labels, each chain has its required order with probability $1/K!$, independently of the other chains. The isolated vertices impose no condition. Equation~\eqref{eq:order-probability} therefore gives the first formula in~\eqref{eq:separated}.

At $j=r+1$, the gamma density gives
\[
\E\left(\frac{G_{K+1}^{\,d}}{K!}\right)^r
=\frac{\Gamma(K+dr+1)}{\Gamma(K+1)(K!)^r}
=\frac{(K+dr)!}{(K!)^{r+1}}.
\]
This is the required count, including its mass at $r=0$. For example, when $K=1$ it is $\Gamma(2+(1-m)r)$, which need not be constant. If $m=K$, every prescription is the same single chain. That separate case gives $N(j)=1$.
\end{proof}

\subsection{Prefixes formed from two extreme blocks}

Part~I~\cite[Proposition~5.28]{BKS1} gives a product
formula when the prefix is $1,2,\ldots,m-1,K$ and
the suffix is increasing. We extend this formula,
for $m\le K/2$, to prefixes formed from two blocks
of the smallest and largest letters, and identify
the corresponding representing measures. These
conditions force successive insertions into nested
intervals. No tree-like assumption is needed.

\begin{theorem}\label{thm:prefix}
Let $2\le m\le K/2$, $1\le a_0<m$, and $b_0=K-m+a_0+1$. Suppose that the actual $m$-prefix of $\pi$ is
\[
1,2,\ldots,a_0,b_0,b_0+1,\ldots,K
\]
and that its actual $m$-suffix $s_1<\cdots<s_m$ is increasing. The letters of the interfix may occur in any order. Set
\[
d=K-m,\qquad \ell=s_{a_0+1}-s_{a_0}-1.
\]
Then
\begin{equation}\label{eq:prefix-count}
N_\pi^m(j)=\prod_{i=0}^{j-1}\binom{id+\ell}{\ell},\qquad
N_\pi^m(r+1)=\E\left(\frac{d^\ell}{\ell!}\prod_{t=1}^{\ell}G_{1+t/d}\right)^r.
\end{equation}
The factors are independent, and the empty product for $\ell=0$ is one.
\end{theorem}
\begin{proof}
We count the choices of values in the first occurrence that can be interleaved with the remaining occurrences. The extreme prefix blocks make the remaining part a smaller packing after standardization. This yields a recurrence, to which~\eqref{eq:binomial-gammas} applies.

Since $2m\le K$, the prefix and suffix occupy disjoint positions. Every suffix value consequently belongs to $\{a_0+1,\ldots,b_0-1\}$. When a new copy is attached, the first $a_0$ shared vertices correspond to its smallest letters, and the other shared vertices to its largest letters. All $d$ new vertices lie strictly between the $a_0$th and $(a_0+1)$st shared vertices. The new suffix lies in the same interval. Repeating this observation proves that every later insertion lies in a nested subinterval.

In the first occurrence, the vertices corresponding to letters at most $s_{a_0}$ precede every later vertex, while those corresponding to letters at least $s_{a_0+1}$ follow every later vertex. The vertices corresponding to the $\ell$ intervening letters form a chain. This chain is incomparable with all $(j-1)d$ vertices added in subsequent copies: it contains no shared vertex, and every comparison with another copy must pass through the first suffix, whose elements corresponding to $s_{a_0}$ and $s_{a_0+1}$ bound the chain.

The later vertices induce a packing of $j-1$ copies with its first prefix removed. In that smaller packing the removed prefix is a forced initial chain of $a_0$ vertices and a forced terminal chain of $m-a_0$ vertices. Removing them does not change the number of linear extensions. An extension of the full poset is therefore determined by an extension of the smaller packing and the $\ell$ places of the intervening first-copy chain among $(j-1)d+\ell$ places. Deleting that chain recovers both choices, so this correspondence is bijective. Thus
\[
N_\pi^m(j)=\binom{(j-1)d+\ell}{\ell}N_\pi^m(j-1),\qquad N_\pi^m(1)=1.
\]
Iteration proves the count. Taking $j=r+1$ in~\eqref{eq:binomial-gammas}, with $h=d$, gives the second formula in~\eqref{eq:prefix-count}. When $\ell=0$, both the interleaving choice and the random product are identically one.
\end{proof}


\section{Rooted trees and moment representations}\label{sec:rooted}

The binomial factors in Proposition~\ref{prop:extreme-endpoint}
count the ways to interleave a branch with the remaining vertices
above its attachment vertex. The same argument applies to rooted
tree posets at larger overlaps, albeit with an additional difficulty that
attaching a new occurrence can subdivide the paths between earlier
attachment vertices. We show that each of these paths eventually
stops changing, and use the resulting count formula to construct
a Stieltjes moment representation for the counting sequence after finitely
many initial terms have been removed. When the branches have
one or two vertices, we also determine the minimum number of
terms that must be removed.

\subsection{Rooted geometry and extension counts}

We call a tree poset \emph{rooted} if it has a unique minimum.
Fix a pattern $\pi\in S_K$ and an admissible overlap $m$ with
\[
K\ge2,\qquad 1\le m\le K/2,
\]
and write $N(j)=N_\pi^m(j)$. We seek a product formula for
these counts whose binomial factors eventually have upper
entries in arithmetic progression. Such a formula will give
the gamma moment representation in Theorem~\ref{thm:rooted-tails}.

By Part~I~\cite[Lemma~3.8 and Remark~5.9]{BKS1}, we may
assume that the prefix and suffix are increasing, since
sorting them preserves the pattern-packing poset up to
order isomorphism. Write them as
\[
P=(p_1<\cdots<p_m),\qquad S=(s_1<\cdots<s_m).
\]
For $j\ge2$, the poset is a tree if and only if
\begin{equation}\label{eq:tree-test}
(p_{i+1}-p_i-1)(s_{i+1}-s_i-1)=0
\qquad(1\le i<m),
\end{equation}
by Part~I~\cite[Theorem~5.23]{BKS1}. It has a unique minimum
if and only if $\min(p_1,s_1)=1$, by
Part~I~\cite[Proposition~5.12]{BKS1}. In the rooted tree cases,
it is a chain precisely when also $\max(p_m,s_m)=K$.
Chains give $N(j)=1$ for all $j$, represented by $\delta_1$.
We therefore assume from now on that, for $j\ge2$, the posets
are rooted trees and are not chains. The corresponding results
for trees with a unique maximum follow by order duality.

Both $p_m$ and $s_m$ are now less than $K$. Choose
\[
L=(\ell_1<\cdots<\ell_m),\qquad
R=(\rho_1<\cdots<\rho_m)
\]
to be $P,S$ in the order for which $\ell_m<\rho_m$.
Interchanging prefix and suffix gives an order isomorphism,
by Part~I~\cite[Lemma~3.10 and Remark~5.9]{BKS1}, so we may
attach occurrences with prefix $L$ and suffix $R$.
This choice does not specify which of $\ell_1,\rho_1$ is $1$.
Number the occurrences from $0$ to $j-1$, and write $(t,a)$
for the vertex corresponding to the pattern letter $a$ in occurrence $t$,
that is, the element in position $a$ on line $t$ of the stack.
The shared vertices satisfy
\[
(t,\rho_i)=(t+1,\ell_i).
\]
Set
\begin{equation}\label{eq:root-parameters}
\begin{gathered}
\alpha_i=\ell_{i+1}-\ell_i-1,\qquad
\beta_i=\rho_{i+1}-\rho_i-1\qquad(1\le i<m),\\
b=K-\rho_m,\qquad c=K-\ell_m,\qquad
\sigma_0=\rho_m-\ell_m=c-b.
\end{gathered}
\end{equation}
Thus $c>b\ge1$ and $\alpha_i\beta_i=0$.

By Part~I~\cite[Definition~5.15, Theorem~5.16 and
Porism~5.17]{BKS1}, the core is a chain, and all branches
go upwards from it and have $b$ vertices.
If $\rho_1=1$, we apply those results
after interchanging prefix and suffix, then transfer the
description through the order isomorphism above. In our chosen
attachment order, the core joins the unique minimum to the
vertex corresponding to $K$ in the last occurrence.

Attaching occurrence $t+1$ makes the vertex
\[
B_t=(t,\rho_m)=(t+1,\ell_m),\qquad 0\le t\le j-2
\]
an attachment vertex.
The $b$ vertices above it in occurrence $t$ form a branch,
while the $c$ vertices above it in occurrence $t+1$ continue
the core. These are all the vertices above the newest
attachment vertex. Since $\min(\ell_1,\rho_1)=1$, any new
vertices below the lowest shared vertex form a chain below
the rest of the poset. The attachment vertices lie in
increasing order on the core, since
\[
B_{t-1}=(t,\ell_m)<(t,\rho_m)=B_t
\qquad(1\le t\le j-2).
\]

An interval between consecutive attachment vertices has
\emph{age} one when its upper vertex becomes an attachment
vertex. Its age increases by one whenever another occurrence
is added. Before occurrence $t+1$ is attached, the segment
from $(t,\ell_m)$ to $(t,\rho_m)$ has $\sigma_0$ cover edges.
Earlier occurrences meet line $t$ only in its prefix, at or
below $(t,\ell_m)$, since a vertex belongs to at most two
adjacent stack lines, by Part~I~\cite[Lemma~5.8]{BKS1}.
They therefore do not subdivide this segment. Each subsequent
attachment uses the same identifications between suffix and
prefix letters, so the length of the interval depends only
on its age. Write $\sigma_a$ for its number of cover edges
at age $a\ge1$.

For $j\ge3$ occurrences, the successive distances read
downwards from the newest attachment vertex are
\[
\sigma_1,\ldots,\sigma_{j-2}.
\]
At $j=2$ there is only one attachment vertex.
Let $C_k$ be the number of vertices strictly above the
$(k+1)$st attachment vertex in this order, excluding its
branch. At the newest attachment vertex, $C_0=c$.
Moving to the preceding one adds the branch of $b$ vertices
and the connecting path, including its upper endpoint. Hence
\begin{equation}\label{eq:root-continuation}
C_k=c+kb+\sum_{a=1}^{k}\sigma_a\qquad(k\ge0).
\end{equation}
This count is independent of $j$, as long as $j\ge k+2$.

The branch at each attachment vertex is incomparable with
these $C_k$ vertices, so it can be interleaved with any
linear extension of them in $\binom{C_k+b}{b}$ ways.
This is the interleaving argument of
Part~I~\cite[Proposition~5.24]{BKS1}, transferred through
the prefix--suffix isomorphism when $\rho_1=1$.
The vertices below the lowest attachment vertex form a chain
below the rest of the poset and contribute no further choice.
Applying the interleaving argument successively gives
\begin{equation}\label{eq:root-count}
N(j)=\prod_{k=0}^{j-2}\binom{C_k+b}{b}\qquad(j\ge1),
\end{equation}
where the empty product gives $N(1)=1$.

The factors in~\eqref{eq:root-count} are therefore determined
by the distances $\sigma_a$. Since
$C_k-C_{k-1}=b+\sigma_k$ for $k\ge1$, their upper entries
will eventually form an arithmetic progression if the
distances $\sigma_a$ eventually stop changing.

Consider an index $i$ with $\alpha_i>0$, which we call
\emph{active}. The tree condition gives $\beta_i=0$, so
$\rho_{i+1}=\rho_i+1$. Before the next occurrence is attached,
the corresponding vertices on the last line form a cover edge.
Indeed, the prefix and suffix are disjoint, so neither vertex
belongs to an earlier line. The new occurrence replaces this
edge by a path with $\alpha_i$ additional vertices.
Inactive indices insert no vertices between their shared
endpoints.

When $B_t$ becomes an attachment vertex, with $t\ge1$,
the path from $B_{t-1}$ to $B_t$ is subdivided at exactly
the active indices in
\[
F_0=\{i:\alpha_i>0,\ \rho_i>\ell_m\}.
\]
After a subdivision at $i$, the inserted path lies between
the letters $\ell_i$ and $\ell_{i+1}$ of the new occurrence.
Its endpoints already belong to two stack lines, so only
its internal vertices can be shared with the following
occurrence. Thus the active indices that can subdivide
this path at the next attachment are exactly those $h$
for which
\[
\ell_i<\rho_h<\rho_{h+1}<\ell_{i+1}.
\]
Call these indices the \emph{children} of $i$, and let
$F_{g+1}$ consist of the children of the indices in $F_g$.
Then $F_g$ records the subdivisions at age $g+1$.
Further changes to the distance can occur only while these
generations remain nonempty.

\begin{lemma}[Stabilization of the distances]\label{lem:rooted-depth}
There are finitely many nonempty generations $F_g$.
If their number is $D$, with $D=0$ when $F_0=\varnothing$, then
\[
D\le\lfloor m/2\rfloor,
\]
and the distances $\sigma_a$ are constant for
$a\ge\max(1,D)$.
\end{lemma}
\begin{proof}
For $1\le h\le m$, let
\[
f(h)=\#\{r\in[m]:\ell_r<\rho_h\}.
\]
This function is nondecreasing. If $h$ is active, then
$\rho_{h+1}=\rho_h+1$, and disjointness of $L$ and $R$ gives
$f(h+1)=f(h)$. An active index $h$ is a child of an active
index $i$ precisely when $f(h)=i$, while $h\in F_0$ precisely
when $h$ is active and $f(h)=m$. Hence every index has at most
one parent, and no index in $F_0$ has an active parent.

Suppose that $f(i)>i$ and that $h$ is a child of $i$.
Monotonicity forces $h<i$. If $h=i-1$, then
$f(i)=f(h+1)=f(h)=i$, a contradiction. Thus $h\le i-2$.
Moreover, $f(h)=i>h$, so this argument applies at every step
of a path starting in $F_0$. Such a path is finite and has
no repeated index. Together with uniqueness of parents,
this shows that the generations form a forest with no
repeated index.

If $D\ge1$, a longest path has $D$ vertices. Its first index
is at most $m-1$, and its indices decrease by at least two
at each step. Therefore
\[
1\le m-1-2(D-1),
\]
which gives the bound. The case $D=0$ is immediate.
Since $F_g$ records the subdivisions at age $g+1$,
there are no further changes once the interval has age
$\max(1,D)$.
\end{proof}

\begin{figure}[htbp]
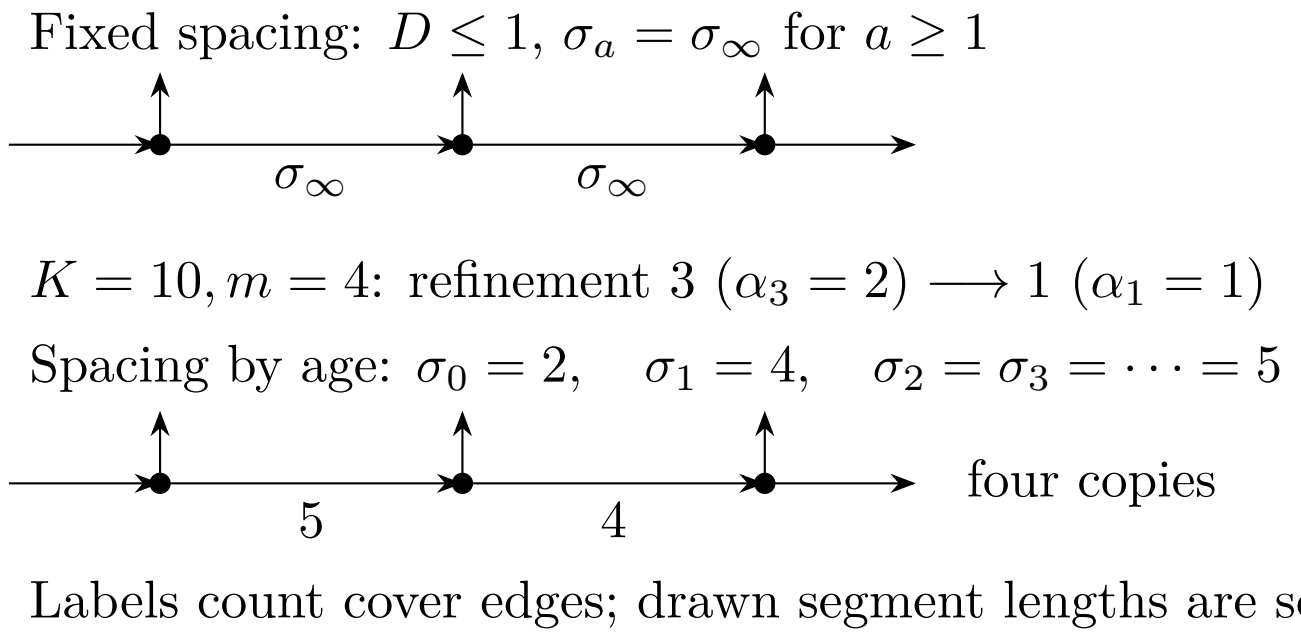

\centering\RootedFigure
\caption{Distances between attachment vertices, measured in
cover edges, with order following the arrows. A branch is
attached at each marked vertex. In the upper example the
distances are fixed. In the lower example, index $3$ inserts
two vertices and its child, index $1$, inserts one more at
the next attachment. An interval therefore has length $4$
at age one and length $5$ at every later age. With four
occurrences, the successive distances from oldest to newest
are $5,4$. The drawing is schematic.}\label{fig:rooted}
\end{figure}

To compute the distances, write
$w_g=\sum_{i\in F_g}\alpha_i$, taking $w_g=0$ when $F_g$
is empty. The total increase at age $g+1$ is $w_g$ cover
edges, so
\begin{equation}\label{eq:spacing}
\sigma_a=\sigma_0+\sum_{g=0}^{a-1}w_g,\qquad
\sigma_\infty=\sigma_0+\sum_{g\ge0}w_g.
\end{equation}
Here $\sigma_\infty$ is the distance after the subdivisions
have stopped, as illustrated in Figure~\ref{fig:rooted}.

We can now evaluate the binomial factors once the distances
have stabilized. To record the slope and intercept of $C_k$, set
\begin{equation}\label{eq:stable-data}
E=\sum_{g\ge0}g w_g,\qquad
\delta=b+\sigma_\infty,\qquad Q=\max(1,D).
\end{equation}
Summing~\eqref{eq:spacing}, for $k\ge Q-1$ we obtain
\[
\sum_{a=1}^{k}\sigma_a
=k\sigma_0+\sum_{g=0}^{k-1}(k-g)w_g
=k\sigma_\infty-E.
\]
When $D\ge1$, the last equality also holds at $k=D-1$, because the
coefficient of $w_{D-1}$ is then zero. Equations
\eqref{eq:root-continuation} and~\eqref{eq:root-count}
therefore give
\begin{align}
C_k&=c-E+k\delta &&(k\ge Q-1),\label{eq:affine}\\
\frac{N(j+1)}{N(j)}
&=\frac1{b!}\prod_{t=1}^{b}
\bigl(c-E+(j-1)\delta+t\bigr)
&&(j\ge Q).\label{eq:root-ratio}
\end{align}
Although the intercept $c-E$ can be negative, every $C_k$
is positive because it counts vertices. The increment
$\delta$ describes these vertex counts and can differ from
$d=K-m$, the difference between the starting positions of
consecutive occurrences.

\subsection{A Stieltjes moment representation}

We use the gamma variables $G_\gamma$ of shape $\gamma$
and rate one introduced in Section~\ref{sec:moments}.

\begin{theorem}\label{thm:rooted-tails}
For a rooted nonchain family as above and an integer $q\ge Q$,
let
\begin{equation}\label{eq:tail-law}
X_q=\frac{\delta^b}{b!}
\prod_{t=1}^{b}G_{(C_{q-1}+t)/\delta},
\end{equation}
where the factors are independent. Then
\[
\frac{N(q+r)}{N(q)}=\E X_q^r\qquad(r\ge0).
\]
Thus $(N(q+r))_{r\ge0}$ is represented by
$N(q)\Law(X_q)$. For every tree family in the stated overlap
range with a unique minimum or a unique maximum, including
chains, the sequence $(N(q+r))_{r\ge0}$ is Stieltjes for some
integer $q$ with $1\le q\le\max(1,\lfloor m/2\rfloor)$.
\end{theorem}
\begin{proof}
Since $q\ge Q$, equation~\eqref{eq:affine} gives
$C_{q-1+u}=C_{q-1}+u\delta$ for every $u\ge0$.
Multiplying the ratios in~\eqref{eq:root-ratio} yields
\begin{align*}
\frac{N(q+r)}{N(q)}
&=\prod_{u=0}^{r-1}\frac1{b!}
\prod_{t=1}^{b}(C_{q-1}+u\delta+t)\\
&=\left(\frac{\delta^b}{b!}\right)^r
\prod_{t=1}^{b}\left(\frac{C_{q-1}+t}{\delta}\right)_r.
\end{align*}
All gamma shapes are positive, since $C_{q-1}\ge c>0$
and $\delta>0$. The gamma moment formula and independence
identify the last expression with $\E X_q^r$.
At $r=0$, both sides are one.

Lemma~\ref{lem:rooted-depth} gives
$Q\le\max(1,\lfloor m/2\rfloor)$. Chains have already
been treated, and order duality preserves the number
of linear extensions, proving the remaining cases.
\end{proof}

The sequence $(N(q+r))_{r\ge0}$ is obtained by deleting
the first $q-1$ terms of the original count sequence,
since $N(q+r)=A_{r+q-1}$.

\subsection{The minimum number of terms to remove}

For branches of length one or two, it turns out that $Q-1$ is the minimum
number of initial terms that must be removed for the
counting sequence to admit a moment representation.

\begin{theorem}\label{thm:exact-start}
For a rooted nonchain family with $b=1$ or $b=2$, and any
integer $q\ge1$,
\[
\begin{aligned}
(N(q+r))_{r\ge0}\text{ is Hamburger}
&\quad\Longleftrightarrow\quad
(N(q+r))_{r\ge0}\text{ is Stieltjes}\\
&\quad\Longleftrightarrow\quad q\ge Q.
\end{aligned}
\]
\end{theorem}
\begin{proof}
\emph{Determinacy for $q\ge Q$.}
Carleman's criterion implies that a Stieltjes moment sequence
$(m_r)_{r\ge0}$ satisfying
\begin{equation}\label{eq:Carleman}
\sum_{r\ge1}m_r^{-1/(2r)}=\infty
\end{equation}
has a unique positive representing measure on $\R$.
For a fixed $q\ge Q$, Theorem~\ref{thm:rooted-tails}
provides Stieltjes existence. Its normalized moments
$m_r=N(q+r)/N(q)$ satisfy, by Stirling's
formula,
\[
\log m_r=br\log r+O(r),\qquad
m_r^{-1/(2r)}\asymp r^{-b/2}.
\]
Thus~\eqref{eq:Carleman} diverges for $b=1,2$.
Multiplying by the positive mass $N(q)$ preserves
determinacy, so every such sequence is Hamburger determinate.

\medskip\noindent
\emph{Excluding $q=D-1$ in the Stieltjes case.}
If $D\le1$, then $Q=1$ and no initial terms need to be
removed. Suppose that $D\ge2$, so $Q=D$, and put
$\Delta=w_{D-1}>0$. Write
\[
\mu_D=N(D)\Law(X_D).
\]
If a positive measure $\nu$ on $[0,\infty)$ represented
$(N(D-1+r))_{r\ge0}$, then $x\nu(dx)$ would
represent $(N(D+r))_{r\ge0}$. Determinacy would
force $x\nu=\mu_D$, and hence
\[
\nu((0,\infty))=\int_{(0,\infty)}x^{-1}\,\mu_D(dx).
\]
Direct integration gives
$\E G_\gamma^{-1}=1/(\gamma-1)$ for $\gamma>1$
and infinity for $0<\gamma\le1$. If any shape in $X_D$
is at most one, Tonelli's theorem and independence make
this inverse moment infinite, contradicting the required
mass $N(D-1)$.

If every shape exceeds one, use
\[
C_{D-1}-C_{D-2}=b+\sigma_{D-1}=\delta-\Delta
\]
to obtain
\begin{align}\label{eq:inverse-obstruction}
\int x^{-1}\,\mu_D(dx)
&=\frac{N(D)b!}
{\prod_{t=1}^{b}(C_{D-1}-\delta+t)}\notag\\
&=\frac{N(D)b!}
{\prod_{t=1}^{b}(C_{D-2}-\Delta+t)}\notag\\
&>\frac{N(D)b!}{\prod_{t=1}^{b}(C_{D-2}+t)}
=N(D-1).
\end{align}
All denominator factors are positive, and subtracting
$\Delta>0$ strictly decreases their product. The resulting
mass on $(0,\infty)$ already exceeds $N(D-1)$, so an atom
at zero cannot supply a representing measure.

A Stieltjes measure for $(N(q+r))_{r\ge0}$ with $q<D-1$
would give the excluded Stieltjes representation of
$(N(D-1+r))_{r\ge0}$ after multiplication by $x^{D-1-q}$.
Thus $(N(q+r))_{r\ge0}$ is not Stieltjes for any $q<D$.

\medskip\noindent
\emph{Excluding $q<D$ in the Hamburger case.}
Suppose that a positive measure $\nu$ on $\R$ represents
$(N(q+r))_{r\ge0}$ for some $q<D$. Choose an even integer $2u$
such that $q+2u\ge D$. The positive measure $x^{2u}\nu(dx)$
represents $(N(q+2u+r))_{r\ge0}$, so Hamburger determinacy gives
\[
x^{2u}\nu(dx)
=N(q+2u)\Law(X_{q+2u})(dx).
\]
The right-hand side is concentrated on $(0,\infty)$.
Since $x^{2u}>0$ on $(-\infty,0)$, the measure $\nu$
has no mass there. It would therefore be a Stieltjes
measure for $(N(q+r))_{r\ge0}$, which we have excluded.
Theorem~\ref{thm:rooted-tails} supplies Stieltjes representations
for all $q\ge Q$.
\end{proof}


\section{Two-level grids and total nonnegativity}\label{sec:grids}

In Sections~\ref{sec:enumeration} and~\ref{sec:rooted}, we
used the tree structure of the posets to count their
linear extensions. In Section~\ref{sec:enumeration},
conditioning on the trunk variables made the remaining
order constraints independent, leading to integral formulas.
For the rooted trees of Section~\ref{sec:rooted}, successive
interleavings reduced to products of binomial coefficients.

\emph{Shifted grid posets}, introduced in Part~I~\cite{BKS1} and recalled below, require a different counting argument. For two-level grids, fixing the values assigned to selected upper-row vertices leaves choices counted by products of binomial coefficients. We organize these choices into a transfer matrix and use it to prove total nonnegativity of the entire Hankel matrix of the packing sequence.

\subsection{From the grid to a transfer matrix}

We recall the shifted-grid construction from Part~I~\cite[Theorem~5.33]{BKS1}. Let $m\ge2$ and suppose that the $m$-prefix and $m$-suffix of $\pi\in S_K$ are the increasing words $p_1\cdots p_m$ and $s_1\cdots s_m$, with
\begin{equation}\label{eq:grid-boundaries}
1=p_1<s_1<p_2<s_2<\cdots<p_m<s_m=K,
\qquad m\le K/2.
\end{equation}
Here the entries are in their order of appearance in $\pi$. Let $d=K-m$ and
\[
h_i=s_i-p_i-1\quad(1\le i\le m),\qquad
v_i=p_{i+1}-s_i-1\quad(1\le i<m).
\]
Thus $h_i$ counts the values strictly between $p_i$ and $s_i$, while $v_i$ counts those between $s_i$ and $p_{i+1}$.

For a packing with $j$ occurrences, represent position $qd+i$ by the vertex
\[
u_{i,q}=(q+i-1,i-1),\qquad 1\le i\le m,\quad0\le q\le j.
\]
In occurrence $q+1$, these vertices appear along the chain in the order
\[
u_{1,q}\longrightarrow u_{1,q+1}\longrightarrow
u_{2,q}\longrightarrow u_{2,q+1}\longrightarrow\cdots
\longrightarrow u_{m,q}\longrightarrow u_{m,q+1}.
\]
Each horizontal arrow $u_{i,q}\to u_{i,q+1}$ contains $h_i$ internal vertices, and each vertical arrow $u_{i,q+1}\to u_{i+1,q}$ contains $v_i$. Placing the other positions of each occurrence along these arrows in increasing order of their values in $\pi$ recovers the chains in~\eqref{eq:poset-chains}. The vertices with $0<q<j$ belong to two consecutive occurrences, whereas those with $q=0$ or $q=j$ belong to only one.

Conversely, any nonnegative subdivision parameters determine such a family. Take
\[
K=2m+\sum_{i=1}^m h_i+\sum_{i=1}^{m-1}v_i,
\qquad p_1=1,
\]
define $s_i=p_i+h_i+1$ and $p_{i+1}=s_i+v_i+1$ successively, and arrange the unused letters in the interfix, which we may take to be increasing. This realizes the shifted grid with the chosen parameters.

For two levels, write
\[
(h,v,k)=(h_1,v_1,h_2)\in\mathbb Z_{\ge0}^3.
\]
In this case $K=h+v+k+4$, $d=K-2$, and the prefix and suffix are $(1,h+v+3)$ and $(h+2,K)$, respectively, as in Part~I~\cite[Corollary~5.36]{BKS1}. Write $x_q=u_{1,q}$ and $y_q=u_{2,q}$. Occurrence $i$ follows the chain from $x_{i-1}$ through $x_i$ and $y_{i-1}$ to $y_i$, with $h$, $v$, and $k$ internal vertices on its three successive parts.

Fix a pattern of this form, and write $N(j)=N_\pi^2(j)$ and $A_r=N(r+1)$. A packing with $j$ occurrences has $n=2+jd$ positions. A linear extension assigns the values $1,\ldots,n$ to the vertices in its total order. Let $a_i$ be the value assigned to $y_{i-1}$, for $1\le i\le j$. Because the upper row is a chain, these values are strictly increasing. Moreover, the final $k$ upper subdivisions and $y_j$ lie above all other vertices, so
\[
a_j=n-k-1.
\]

\begin{figure}[htbp]
\centering
\begin{tikzpicture}[x=10mm,y=10mm,>=Stealth,font=\small,
  corner/.style={circle,draw,fill=white,inner sep=2.1pt},
  ordinary/.style={->,line width=.65pt,draw=black!55},
  lowerprefix/.style={->,line width=1.65pt,draw=black},
  verticalchoice/.style={->,line width=1.55pt,draw=blue!65!black},
  upperchoice/.style={->,line width=1.55pt,draw=orange!70!black,dashed}]

\foreach \q/\X in {0/0,1/2.15,2/4.30,3/6.45}{
  \node[corner] (x\q) at (\X,0) {};
  \node[below=4pt] at (x\q) {$x_{\q}$};
}
\foreach \q/\X in {0/2.15,1/4.30,2/6.45,3/8.60}{
  \node[corner] (y\q) at (\X,1.75) {};
  \node[above=3pt] at (y\q) {$y_{\q}$};
}
\node[above=21pt] at (y0) {$a_1$};
\node[above=21pt] at (y1) {$a_2$};
\node[above=21pt] at (y2) {$a_3$};

\draw[lowerprefix] (x0) -- node[above=4pt] {$h$} (x1);
\draw[lowerprefix] (x1) -- node[above=4pt] {$h$} (x2);
\draw[ordinary] (x2) -- node[above=4pt] {$h$} (x3);

\draw[ordinary] (x1) -- node[left] {$v$} (y0);
\draw[verticalchoice] (x2) -- node[left] {$v$} (y1);
\draw[ordinary] (x3) -- node[left] {$v$} (y2);

\draw[ordinary] (y0) -- node[below] {$k$} (y1);
\draw[upperchoice] (y1) -- node[below] {$k$} (y2);
\draw[ordinary] (y2) -- node[below] {$k$} (y3);

\draw[line width=.8pt] (0,-.76)--(0,-.62);
\draw[line width=.8pt] (0,-.69)--(4.30,-.69);
\draw[line width=.8pt] (4.30,-.76)--(4.30,-.62);
\node[below=1pt,font=\footnotesize] at (2.15,-.76)
  {$2(h+1)+1$ lower-row vertices};
\end{tikzpicture}
\caption{The first three occurrences of a two-level shifted grid.
The numbers $a_1,a_2,a_3$ are the values assigned to $y_0,y_1,y_2$.
At the second step, the highlighted lower prefix precedes the
vertical subdivision $x_2\to y_1$, while the dashed upper segment
$y_1\to y_2$ has its internal vertices assigned values between $a_2$ and $a_3$.
Each edge label gives its number of internal subdivision vertices.}
\label{fig:grid}
\end{figure}
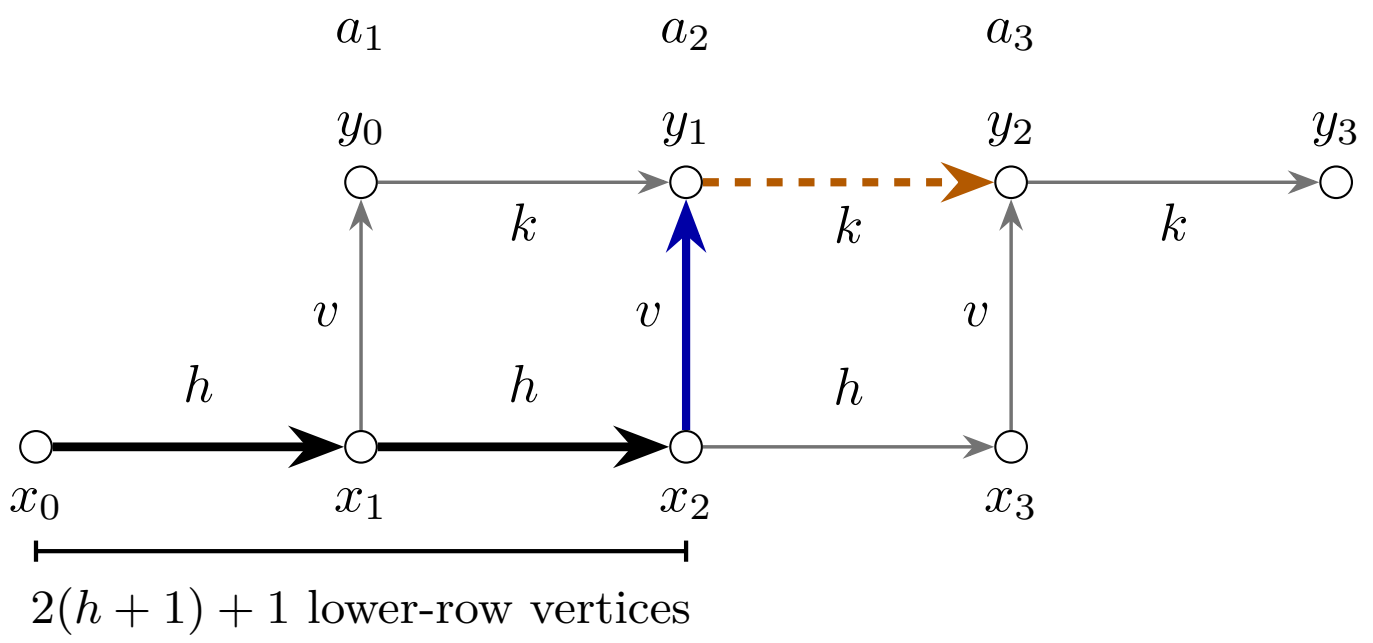

We count linear extensions with $a_1,\ldots,a_j$ fixed.
Each upper-row edge $y_{i-1}\to y_i$ has $k$ internal vertices,
whose assigned values can be chosen in
\[
\binom{a_{i+1}-a_i-1}{k}
\]
ways. These choices are independent over $1\le i<j$, since the
intervals between values assigned to successive upper-row vertices are disjoint.

After choosing the upper subdivisions, we assign values to the
vertical subdivisions $x_i\to y_{i-1}$, in order of increasing $i$.
At step $i$, the earlier upper-row vertices, upper subdivisions and
vertical subdivisions occupy $(i-1)(k+1+v)$ values smaller than $a_i$.
There are thus $a_i-1-(i-1)(k+1+v)$ values smaller than $a_i$ still available.
The lower chain through $x_i$ contains $i(h+1)+1$ vertices, all of
which must precede the $v$ internal vertices on $x_i\to y_{i-1}$.
The smallest $i(h+1)+1$ still-available values less than $a_i$ are thus
reserved for the lower row, as shown for $i=2$ in
Figure~\ref{fig:grid}. We may choose any $v$ of the other available
values smaller than $a_i$ for the current vertical subdivision, giving
\[
\binom{a_i-1-(i-1)(k+1+v)-i(h+1)-1}{v}
=\binom{a_i-d(i-1)-h-3}{v}.
\]
The lower prefix required at step $i$ is contained in the prefix
required at step $i+1$, so these choices are compatible.
After the terminal upper chain is placed, the unused values belong
to the lower row and are assigned in increasing order.
Conversely, each linear extension specifies the same upper and
vertical choices of values.

With the convention $\binom{b}{c}=0$ whenever $b<c$ and $c\ge0$, so that impossible choices of values contribute zero, we obtain
\begin{equation}\label{eq:rank-sum}
\begin{aligned}
N(j)&=\sum_{1\le a_1<\cdots<a_j=n-k-1}
\prod_{i=1}^j\binom{a_i-d(i-1)-h-3}{v}\times
\prod_{i=1}^{j-1}\binom{a_{i+1}-a_i-1}{k}.
\end{aligned}
\end{equation}

The dependence on $i$ in the first binomial factor disappears after the change of variables
\[
t_i=a_i-d(i-1)-h-v-2.
\]
A nonzero term has $t_i\ge1$, and $a_j=n-k-1$ becomes $t_j=1$. Define
\begin{equation}\label{eq:grid-transfer}
\begin{aligned}
\Delta&=h+v+1=d-k-1,&\qquad w(t)&=\binom{t+v-1}{v},\\
T_{s,t}&=\binom{d+t-s-1}{k}w(t),&&s,t\ge1.
\end{aligned}
\end{equation}
The contribution of $t_1,\ldots,t_j$ in~\eqref{eq:rank-sum} is then $w(t_1)T_{t_1,t_2}\cdots T_{t_{j-1},t_j}$. Since $t_j=1$, summing over the other states gives
\begin{equation}\label{eq:source}
N(j)=wT^{j-1}e_1,\qquad A_r=wT^re_1,
\end{equation}
where $w$ is a row vector and $e_1$ is the first coordinate column. The binomial convention also enforces the strict increase of the original values $a_i$, since a nonzero transition requires $a_{i+1}-a_i\ge k+1$.

Although $T$ is infinite, each expression in~\eqref{eq:source} is a finite sum. A transition from state $s$ to state $t$ has nonzero weight only if $t\ge s-\Delta$. Working backward from the final state $t_j=1$, we find
\begin{equation}\label{eq:finite-states}
1\le t_i\le1+(j-i)\Delta.
\end{equation}
There is also a distinction between the first state and subsequent ones: the first contributes $w(t_1)$ without an incoming $T$-transition. To put all steps into one matrix, adjoin state 0 and set
\[
S_{0,t}=w(t)\quad(t\ge1),\qquad
S_{s,t}=T_{s,t}\quad(s,t\ge1),\qquad
S_{s,0}=0\quad(s\ge0).
\]
Then $N(j)=(S^j)_{0,1}$. If $k=0$, the first row of $T$ already equals $w$, and the formula simplifies to $N(j)=(T^j)_{1,1}$. For $k>0$, the source row must be treated separately.

\subsection{From the transfer matrix to Hankel positivity}

A matrix is \emph{totally nonnegative} if all its finite square minors are nonnegative, with rows and columns in their prescribed increasing order. For a real sequence $a=(a_r)_{r\ge0}$, write
\[
H_n(a)=(a_{i+j})_{0\le i,j\le n},\qquad
H_n^{(1)}(a)=(a_{i+j+1})_{0\le i,j\le n}.
\]
The Hamburger criterion characterizes moment sequences on $\mathbb R$ by positive semidefiniteness of all the $H_n(a)$. For a representing measure on $[0,\infty)$, the Stieltjes criterion requires the same condition for both $H_n(a)$ and $H_n^{(1)}(a)$. Total nonnegativity of the full Hankel matrix $(a_{p+q})_{p,q\ge0}$ implies both conditions, because the shifted Hankel matrix is obtained by omitting its first row.

\begin{theorem}\label{thm:grid}
For every two-level shifted-grid family with $h,v,k\ge0$, the infinite Hankel matrix $(A_{p+q})_{p,q\ge0}$ is totally nonnegative. In particular, $(A_r)_{r\ge0}$ is a Stieltjes probability moment sequence.
\end{theorem}

\begin{proof}
We begin with the binomial kernel in~\eqref{eq:grid-transfer}. On integer indices, let $U_{a,b}=\ind_{\{b\ge a\}}$. In every finite submatrix the rows of $U$ are nested suffixes of ones. Subtracting each selected row from its predecessor leaves disjoint intervals of ones, in increasing order, so every square minor is zero or one. Thus $U$ is totally nonnegative. Its powers have entries
\[
(U^{k+1})_{a,b}=\binom{b-a+k}{k}\quad(b\ge a),
\]
and vanish for $b<a$. Indeed, matrix multiplication counts weakly increasing intermediate indices between $a$ and $b$. They lie in a finite interval, so the finite Cauchy--Binet formula also shows that every power of $U$ is totally nonnegative.

Since $\Delta=d-k-1$, the unweighted kernel of $T$ is
\[
C_{s,t}:=\binom{d+t-s-1}{k}=(U^{k+1})_{s,t+\Delta}.
\]
The shift of column indices preserves their order and hence the signs of all minors. Multiplying column $t$ by $w(t)>0$ gives $T$, proving that $T$ is totally nonnegative.

We need the same property after adjoining the source row $w$. Extend the row index of $C$ to all integers. For any fixed finite set of columns, dividing row $-L$ by $L^k/k!$ and letting $L\to\infty$ gives the constant row $(1,1,\ldots)$. When $k=0$, the divisor is one. Because row $-L$ precedes every positive-indexed row, all minors involving this row are nonnegative before taking the limit and remain so afterwards. We may therefore adjoin the constant row above $C$. Scaling its columns by $w(t)$ produces the source row $w$, and adjoining the zero column indexed by 0 yields $S$. Thus $S$ is totally nonnegative as well.

The passage from a transfer matrix to the rows generated by its powers is a standard production-matrix argument~\cite[Proposition~2.8 and Theorem~2.9]{Sokal2023}. To see it in finite dimensions, let $P$ be a totally nonnegative square matrix and $e$ its first coordinate column. If $O_{n-1}$ consists of the first $n$ rows of the matrix with rows $e^{\mathsf T}P^q$, $q\ge0$, then
\begin{equation}\label{eq:output-factorization}
\begin{pmatrix}e^{\mathsf T}\\O_{n-1}P\end{pmatrix}
=\begin{pmatrix}1&0\\0&O_{n-1}\end{pmatrix}
\begin{pmatrix}e^{\mathsf T}\\P\end{pmatrix}.
\end{equation}
The first factor is totally nonnegative by induction, beginning with the single row $e^{\mathsf T}$. A minor of the second factor that avoids its first row is a minor of $P$. If it includes that row, it vanishes unless the first column is also selected, in which case expansion along the row leaves a minor of $P$ with positive sign. Cauchy--Binet proves that the full matrix of output rows is totally nonnegative. Applying the same argument to $P^{\mathsf T}$ and transposing gives the corresponding assertion for the columns $P^qe$.

Apply this result to principal finite truncations of $S$ and $T$. For any fixed terminal state and number of transitions, the bound on the decrease at each transition gives a bound on all preceding states. This includes paths beginning at source 0, since the first positive state must be able to reach the terminal state in the remaining steps. Consequently, each relevant entry stabilizes once the truncation is large enough. Passing to these values shows that the matrices
\[
B_{p,t}=(wT^p)_t,\qquad E_{t,q}=(T^q)_{t,1},
\qquad p,q\ge0,\quad t\ge1,
\]
are totally nonnegative. Here $B$ is obtained from the row-output matrix of $S$ by omitting its first row and column 0, while $E$ is the column-output matrix of $T$.

They factor the Hankel matrix, because~\eqref{eq:source} gives
\[
(BE)_{p,q}=wT^{p+q}e_1=A_{p+q}.
\]
The products involved in this factorization are finite. Indeed, $E_{t,q}=0$ for $t>1+q\Delta$, so a fixed finite selection of columns of $E$ admits a common bound on the intermediate indices. Every minor of $BE$ is therefore a finite Cauchy--Binet sum of products of nonnegative minors of $B$ and $E$. This proves total nonnegativity of the Hankel matrix.

In particular, every principal minor of $H_n(A)$ and $H_n^{(1)}(A)$ is nonnegative. Both matrices are symmetric and hence positive semidefinite, so the Stieltjes criterion gives a representing measure on $[0,\infty)$. Its total mass is $A_0=w(1)=1$.
\end{proof}

\subsection{Catalan and Fuss--Catalan laws}

Theorem~\ref{thm:grid} gives existence of a representing measure without generally identifying it. For some of these grids, however, the packing counts are known explicitly. For integers $p\ge1$ and $j\ge0$, define the Fuss--Catalan numbers by
\[
F_j^{(p)}=\frac1{pj+1}\binom{(p+1)j}{j}.
\]
The Catalan numbers are $F_j^{(1)}$.

For $K\ge4$, consider a pattern with prefix $13$, suffix $2K$, and intervening entries in any order. Its two-level grid has parameters $(h,v,k)=(0,0,K-4)$, and Part~I~\cite[Proposition~5.37]{BKS1} gives, at overlap 2,
\[
N_\pi^2(j)=F_j^{(K-3)}.
\]
This family is covered by Theorem~\ref{thm:grid}, so the point here is to identify its representing measure.

The same Fuss--Catalan counts occur for two further classes, whose enumerations do not follow from the two-level theorem. For $a\ge2$, the patterns
\[
1(a+1)2(a+2)\cdots a(2a),\qquad
(a+1)1(a+2)2\cdots(2a)a
\]
both have counts $F_j^{(1)}$ at overlap $m=2a-2$,
as established in
\cite[Proposition~2]{BlitvicKammounSteingrimsson2023}. Since $d=K-m=2$, every prescribed occurrence begins at an odd position. Exchanging the entries in each consecutive pair of positions therefore transforms packings of the first pattern into packings of the second, and conversely.

For $m\ge2$, the pattern
\[
\pi=1\,3\,4\cdots(m+2)\,2\,(m+3)\cdots(2m+1)
\]
has, at overlap $m$, the enumeration $N_\pi^m(j)=F_j^{(m)}$. This enumeration was obtained through dense clusters by
Elizalde and Noy~\cite[Theorem~4.2]{ElizaldeNoy2012}
and also appears in
\cite[Proposition~3]{BlitvicKammounSteingrimsson2023}.
As explained in Part~I~\cite[discussion preceding
Proposition~5.37]{BKS1}, the corresponding packing
poset differs from the two-level grid poset for
$1\,3\,4\cdots(m+2)\,2\,(m+3)$ at overlap $2$
by a terminal chain of $m-2$ vertices, which does
not change its number of linear extensions. The restriction $m\ge2$ matters. For $m=1$ the pattern is $132$, whose two-occurrence count is 3 rather than $F_2^{(1)}=2$.

For each family the packing moments are $A_r=F_{r+1}^{(p)}$,
with the appropriate value of $p$. The measures with moments
$F_r^{(p)}$ and $F_{r+1}^{(p)}$ are, respectively, the
Fuss--Catalan distribution of order $p$ and the diagonal
Raney distribution with parameters $(p+1,p+1)$, studied
in~\cite{PensonZyczkowski2011}. We use the beta-product
realization implicit in the Mellin-convolution factorization
of that paper (Appendix~B). Recall that $\operatorname{Beta}(\alpha,\beta)$, for $\alpha,\beta>0$, is the probability law on $(0,1)$ with density
\[
\frac{x^{\alpha-1}(1-x)^{\beta-1}}{\mathrm B(\alpha,\beta)},
\qquad
\mathrm B(\alpha,\beta)=\frac{\Gamma(\alpha)\Gamma(\beta)}{\Gamma(\alpha+\beta)}.
\]
Its $r$th moment is $(\alpha)_r/(\alpha+\beta)_r$. For each integer $p\ge1$, take independent random variables
\[
B_i\sim\operatorname{Beta}\left(
\frac{i}{p+1},\frac{i+1}{p}-\frac{i}{p+1}\right),
\qquad 1\le i\le p,
\]
and set
\[
X_p=\frac{(p+1)^{p+1}}{p^p}\prod_{i=1}^p B_i,
\qquad \mu_p=\Law(X_p).
\]
Both parameters of each beta law are positive.

\begin{corollary}\label{cor:fuss}
For each of the three families above, with its indicated value of $p$, the packing moments $A_r=F_{r+1}^{(p)}$ are represented by the probability measure $x\mu_p(dx)$. For $p=1$, this packing law has density
\[
\frac1{2\pi}\sqrt{x(4-x)}\,\ind_{(0,4)}(x).
\]
\end{corollary}

\begin{proof}
The beta moment formula and independence give
\begin{align*}
\E X_p^r
&=\left(\frac{(p+1)^{p+1}}{p^p}\right)^r
\prod_{i=1}^p\frac{(i/(p+1))_r}{((i+1)/p)_r}\\
&=\frac{\Gamma((p+1)r+1)}{\Gamma(r+1)\Gamma(pr+2)}
=F_r^{(p)}.
\end{align*}
For the second equality, the products of rising factorials satisfy
\[
\prod_{i=1}^{p+1}\left(\frac{i}{p+1}\right)_r
=\frac{((p+1)r)!}{(p+1)^{(p+1)r}},
\qquad
\prod_{i=1}^p\left(\frac{i+1}{p}\right)_r
=\frac{(pr+1)!}{p^{pr}}.
\]
Removing the numerator factor $(1)_r=r!$ and multiplying by the displayed scale factor proves the claimed moment identity.

Since $\int x\,\mu_p(dx)=F_1^{(p)}=1$, weighting $\mu_p$ by $x$ gives a probability measure whose $r$th moment is $F_{r+1}^{(p)}$. When $p=1$, we have $X_1=4B_1$ with $B_1\sim\operatorname{Beta}(1/2,3/2)$. The density of $X_1$ is $\sqrt{(4-x)/x}/(2\pi)$ on $(0,4)$, and multiplication by $x$ gives the asserted packing density.
\end{proof}

The unweighted measures $\mu_p$ are the Fuss--Catalan distributions studied in~\cite{PensonZyczkowski2011}.


\section{General tree posets}\label{sec:nonrooted}

We now revisit the tree posets in the case where the first and last letters are not consecutive, that is, $|\pi_1-\pi_K|>1$. Throughout, we consider only patterns $\pi\in S_K$ with an admissible overlap
$m\le K/2$ satisfying
the tree criterion~\eqref{eq:tree-test}.

\subsection{The trunk integral and its normalization}

\begin{figure}[htbp]
\centering
\begin{tikzpicture}[x=13.5mm,y=5.8mm,
  vertex/.style={circle,draw,fill=white,
    inner sep=1.7pt,font=\footnotesize},
  trunk/.style={vertex,very thick,fill=black!8},
  shared/.style={trunk,fill=black!25},
  internal/.style={trunk,fill=orange!55}]

\node[vertex] (p2) at (-1,0) {$2$};
\node[vertex] (p1) at (-1,1) {$1$};
\node[vertex] (p3) at (-1,2) {$3$};
\node[vertex] (p4) at (-1,4) {$4$};

\node[vertex] (p6) at (1,2) {$6$};
\node[vertex] (p8) at (-1,6) {$8$};

\node[vertex] (p10) at (1,4) {$10$};
\node[vertex] (p11) at (1,6) {$11$};
\node[vertex] (p13) at (1,7) {$13$};
\node[vertex] (p12) at (1,8) {$12$};

\node[shared]   (p5) at (0,3) {$5$};
\node[internal] (p7) at (0,4) {$7$};
\node[shared]   (p9) at (0,5) {$9$};

\draw[very thick] (p5)--(p7)--(p9);

\draw (p2)--(p1)--(p3)--(p5)--(p4);
\draw (p6)--(p5) (p9)--(p8);
\draw (p10)--(p9)--(p11)--(p13)--(p12);

\end{tikzpicture}

\caption{The Hasse diagram for three overlap-one
occurrences of $21354$, with order increasing upward.
The highlighted trunk is $5<_P7<_P9$.
Consecutive occurrences meet at the grey vertices
$5$ and $9$, while the intervening trunk vertex $7$
belongs only to the middle occurrence.}
\label{fig:21354-trunk}
\end{figure}
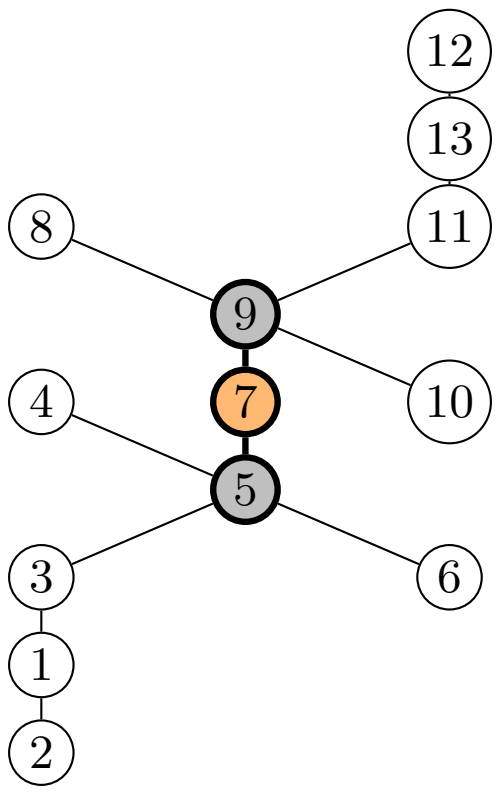

In Section~\ref{sec:enumeration}, we obtained the
adjacent-endpoint formula by first integrating the
chains attached to the trunk and then using symmetry
to reduce the remaining integral to two variables
(Proposition~\ref{prop:integral}).
The first step extends to the present setting.
When the first and last letters are not consecutive, however,
the trunk may contain vertices belonging to only
one occurrence, and the symmetry reduction need
not apply (see Figure~\ref{fig:21354-trunk}).

Concretely, tree posets have the following general structure.
Let $p_1<\cdots<p_m$ and $s_1<\cdots<s_m$
be the sorted letters of the prefix and suffix of $\pi$, and let
\[
\begin{aligned}
a&=\min(p_1,s_1)-1,
&\bar a&=\max(p_1,s_1)-1,\\
b&=K-\max(p_m,s_m),
&\bar b&=K-\min(p_m,s_m).
\end{aligned}
\]

For $r+1$ occurrences, $r\ge1$, the path connecting
all vertices common to successive occurrences is a
chain, which we call the trunk. Here it runs only between the
first and last shared vertices, whereas the conditioning path in Section~\ref{sec:enumeration}
also included the first and last endpoint positions.
Below and above the trunk are chains of lengths
$\bar a$ and $\bar b$, respectively. At each of the
$r$ joins, a lower chain of length $a$ is attached
to the least common vertex, and an upper chain
of length $b$ to the greatest. See Figure~\ref{fig:21354-trunk} for an illustration.

This decomposition follows from the branch description in Part~I \cite[Theorem~5.16 and Porism~5.17]{BKS1} (after the symmetry reductions of \cite[Section~3]{BKS1}). The core is a lattice by \cite[Theorem~5.18]{BKS1}, and therefore a chain when the Hasse diagram is a tree. The trunk is a subchain of this core, excluding its initial and final portions, which have the same form as branches.  Subtracting the attached chains from the total $K+rd$ vertices gives the trunk length
\begin{equation}\label{eq:trunk-length}
L_{r+1}
=K+rd-\bar a-\bar b-r(a+b).
\end{equation}

Index the trunk vertices in increasing order by
$1,\ldots,L_{r+1}$. At the $t$-th join,
$1\le t\le r$, let $i_t^-$ and $i_t^+$ denote
the indices of the least and greatest common
vertices, respectively. Let $x_i$ be the coordinate
of the $i$-th trunk vertex, so that
\[
0<x_1<\cdots<x_{L_{r+1}}<1.
\]

For example, in Figure~\ref{fig:21354-trunk},
we have $a=b=1$ and $\bar a=\bar b=3$.
The trunk consists of the vertices $5<7<9$,
so $L_3=3$. The two joins occur at its first
and third vertices, giving
$i_1^-=i_1^+=1$ and $i_2^-=i_2^+=3$.

The following integral formula extends
Proposition~\ref{prop:integral}. 

\begin{proposition}\label{prop:trunk-integral}
For $r\ge1$, the packing counts
$A_r=N_\pi^m(r+1)$ satisfy
\begin{equation}\label{eq:trunk-integral}
\begin{aligned}
A_r
=\frac{(K+dr)!}
{\bar a!\bar b!(a!b!)^r}\,\,
\int_{0<x_1<\cdots<x_{L_{r+1}}<1}
x_1^{\bar a}(1-x_{L_{r+1}})^{\bar b}
\prod_{t=1}^r
x_{i_t^-}^a(1-x_{i_t^+})^b
\,dx_1\cdots dx_{L_{r+1}}.
\end{aligned}
\end{equation}
\end{proposition}

\begin{proof}
By~\eqref{eq:order-probability},
$A_r=(K+dr)!\Prob(E_P)$.
Condition on the ordered coordinates of the trunk.
As in Proposition~\ref{prop:integral}, an ordered
chain of $k$ vertices below a coordinate $x$
contributes $x^k/k!$, while one above $x$
contributes $(1-x)^k/k!$.

Reasoning as before, the chains below and above the trunk contribute
\[
\frac{x_1^{\bar a}(1-x_{L_{r+1}})^{\bar b}}
{\bar a!\bar b!},
\]
and those attached at the $t$-th join contribute
\[
\frac{x_{i_t^-}^a(1-x_{i_t^+})^b}{a!b!}.
\]
These chains have disjoint vertex sets, so their
conditional probabilities multiply. Integrating
over the ordered trunk coordinates gives
\eqref{eq:trunk-integral}.
\end{proof}

\begin{remark}
The same argument also applies when $m=1$ and
$|\pi_1-\pi_K|=1$. The trunk then consists only of
the $r$ vertices at the joins, and
$\bar a=a+1$, $\bar b=b+1$.
In this case, \eqref{eq:trunk-integral} agrees with
Proposition~\ref{prop:integral} after integrating
the latter's two endpoint coordinates.
\end{remark}

\subsection{A representing measure when $|\pi_1-\pi_K|=2$}

The integral in Proposition~\ref{prop:trunk-integral}
gives each packing count, but is not manifestly
a moment integral. We can nevertheless exhibit
a Stieltjes representing measure for a class of
overlap-one trees whose first and last letters differ by two.

We call a tree family \emph{fully non-rooted} if
its packing posets with at least two occurrences
have neither a unique minimal nor a unique maximal
element. By~\cite[Proposition~5.12]{BKS1},
this means that both $1$ and $K$ belong to
the interfix. At overlap one,
the condition reduces to
\[
1<\pi_1,\pi_K<K.
\]
Therefore, these are precisely the overlap-one families
having neither a unique minimum nor a unique maximum.
Section~\ref{sec:rooted} treats the cases with either unique extremum.

Suppose now that the family is fully non-rooted, $m=1$, and $|\pi_K-\pi_1|=2$.
After reversing the pattern if necessary, take
$\pi_K=\pi_1+2$. With the notation introduced
above,
\[
a=\pi_1-1\ge1,\qquad b=K-\pi_K\ge1.
\]
Consequently, $K=a+b+3$, $d=a+b+2$,
$\bar a=a+2$, and $\bar b=b+2$.

In this case, there being a single letter
between $\pi_1$ and $\pi_K$ greatly simplifies the calculus. For $r+1$
occurrences, the trunk has $r$ vertices shared
by successive occurrences, with one additional
vertex between each pair of successive joins.
Write $0<t_1<\cdots<t_r<1$ for the coordinates
at the joins. The additional vertex between
$t_i$ and $t_{i+1}$ can lie anywhere in that
interval, so integrating its coordinate contributes
$t_{i+1}-t_i$.

The chains attached at each join contribute
$t_i^a(1-t_i)^b$, while the terminal chains
contribute $t_1^{a+2}(1-t_r)^{b+2}$.
Thus, writing $w(t)=t^a(1-t)^b$,
Proposition~\ref{prop:trunk-integral} gives
\begin{equation}\label{eq:grouped-integral}
I_r=
\int_{0<t_1<\cdots<t_r<1}
t_1^{a+2}(1-t_r)^{b+2}
\prod_{i=1}^r w(t_i)
\prod_{i=1}^{r-1}(t_{i+1}-t_i)
\,dt_1\cdots dt_r,
\qquad r\ge1,
\end{equation}
where the packing counts satisfy
\begin{equation}\label{eq:gap-count}
A_r=
\frac{(K+dr)!}
{(a+2)!(b+2)!(a!b!)^r}\,I_r.
\end{equation}

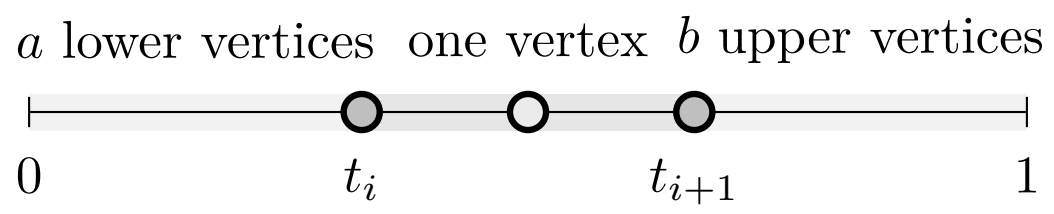
\begin{figure}[htbp]
\centering
\begin{tikzpicture}[
  x=10mm,y=10mm,
  vertex/.style={circle,draw,fill=white,inner sep=2.4pt},
  trunk/.style={vertex,very thick,fill=black!8},
  shared/.style={trunk,fill=black!25},
  every node/.style={font=\small}
]

\fill[black!5] (0,-0.12) rectangle (2.2,0.12);
\fill[black!10] (2.2,-0.12) rectangle (4.4,0.12);
\fill[black!5] (4.4,-0.12) rectangle (6.6,0.12);

\draw (0,0) -- (6.6,0);
\draw (0,-0.1) -- (0,0.1);
\draw (6.6,-0.1) -- (6.6,0.1);

\node[shared] at (2.2,0) {};
\node[trunk] at (3.3,0) {};
\node[shared] at (4.4,0) {};

\node at (1.1,0.48) {$a$ lower vertices};
\node at (3.3,0.48) {one vertex};
\node at (5.5,0.48) {$b$ upper vertices};

\node[below=5pt] at (0,0) {$0$};
\node[below=5pt] at (2.2,0) {$t_i$};
\node[below=5pt] at (4.4,0) {$t_{i+1}$};
\node[below=5pt] at (6.6,0) {$1$};

\end{tikzpicture}
\caption{The allowed coordinate intervals for an occurrence
between consecutive joins $t_i<t_{i+1}$. There are $a$
vertices below $t_i$, one trunk vertex between
$t_i$ and $t_{i+1}$, and $b$ vertices above
$t_{i+1}$. Integrating the intervening trunk
coordinate contributes $t_{i+1}-t_i$.}
\label{fig:gap-endpoints}
\end{figure}

Figure~\ref{fig:gap-endpoints} illustrates the
coordinate intervals underlying
\eqref{eq:grouped-integral}.

Our remaining task is to turn~\eqref{eq:gap-count} into a moment
representation. We do this by separating the integral from the
factorial terms. For the integral contribution, define
\begin{equation}\label{eq:gap-bstar}
b_0^*=1,\qquad
b_r^*=\frac{K!(2r+2)!}{2(a+2)!(b+2)!}\,I_r
\qquad(r\ge1).
\end{equation}
Using~\eqref{eq:gap-count} for $r\ge1$,
and noting that $A_0=N_\pi^m(1)=b_0^*=1$,
we obtain
\begin{equation}\label{eq:gap-factorization}
A_r=\frac{2(K+dr)!}{K!(a!b!)^r(2r+2)!}\,b_r^*,
\qquad r\ge0.
\end{equation}
We will show that $(b_r^*)_{r\ge0}$ consists of the even
moments of a compactly supported symmetric variable $Z$.
The factorial ratio will be identified as the moments of an
independent positive variable $Y$, so that
$A_r=\E(YZ^2)^r$.

The construction of $Z$ rests on the following elementary
lemma about a second-order differential equation.

\begin{lemma}\label{lem:increasing-speed}
Let $c$ be continuous, nonnegative and nondecreasing on $[0,T]$.
For $0\le s\le x\le T$, let $h(x,s;t)$ solve
\[
h_{xx}+t^2c(x)^2h=0,\qquad
h(s,s;t)=1,\qquad h_x(s,s;t)=0.
\]
Then $h(x,s;\cdot)$ is the Fourier transform of a symmetric
probability measure supported on
\[
\left[-\int_s^x c(v)\,dv,\ \int_s^x c(v)\,dv\right].
\]
These measures depend weakly continuously on $(x,s)$.
\end{lemma}

\begin{proof}
Suppose first that $c$ is a positive nondecreasing step
function. A particle starts at the origin at time $s$ with
velocity $c(s)$ or $-c(s)$, each with probability $1/2$.
It moves at constant velocity between changes of speed.
When the speed increases from $c_-$ to $c_+$, it retains
its previous direction with probability $c_-/c_+$, but otherwise
it chooses a new direction uniformly. If $V_-$ and $V_+$
are its velocities before and after the change, then
\[
\E[V_+\mid V_-]=V_-.
\]

Let $X_x$ be its position and $f(x,t)=\E e^{itX_x}$.
Between speed changes,
\[
f_x(x,t)=it\,\E[V_xe^{itX_x}],\qquad
f_{xx}(x,t)=-t^2c(x)^2f(x,t).
\]
At each change, $X_x$ remains fixed, and the conditional
expectation above shows that $f_x$ is continuous. Thus $f$
solves the stated differential equation with $f(s,t)=1$
and $f_x(s,t)=0$, giving $f=h$. Reversing every velocity
leaves the construction unchanged in law, so the position
is symmetric. Its absolute displacement is bounded by
$\int_s^x c(v)\,dv$.

For general $c$, choose positive nondecreasing step functions
converging uniformly to $c$. The corresponding solutions
converge to $h$ by continuous dependence for the differential
equation. Their probability measures have support in a common
compact interval, so subsequential weak limits exist and have
Fourier transform $h$ and the claimed support. Uniqueness of
Fourier transforms identifies the limit. The same compactness
argument, together with the continuity of $h(x,s;t)$ in
$(x,s)$, gives weak continuity of the measures.
\end{proof}

\begin{theorem}\label{thm:gap-two}
For every fully non-rooted overlap-one family with
$|\pi_1-\pi_K|=2$, the packing counts have the representation
\[
A_r=\E(YZ^2)^r,\qquad r\ge0,
\]
where $Y$ and $Z$ are independent random variables. Specifically, 
\begin{equation}\label{eq:gap-Y}
Y=\frac{d^d}{4a!b!}\,U
\prod_{\substack{1\le t\le d\\t\notin\{1,d-1\}}}
G_{1+(t+1)/d},
\end{equation}
where
\[
U\sim\operatorname{Beta}(1+2/d,1/2-2/d)
\quad\text{if }d>4,\qquad U=1\quad\text{if }d=4
\]
and the $G_\alpha$ are rate-one gamma variables of shape
$\alpha$. (All factors are independent.) Furthermore, the variable $Z$ can be chosen symmetric, supported on
$[-L,L]$, where
\[
L=\int_0^1 x^{a/2}(1-x)^{b/2}\,dx,
\]
and with even moments $\E Z^{2r}=b_r^*$ for every $r\ge0$.
\end{theorem}

\begin{proof}
We first express the integrals defining $b_r^*$ through a generating
function. Set
\[
\begin{aligned}
u_0(x)=\frac{x^{a+2}}{(a+1)(a+2)},\qquad
u_{r+1}(x)=\int_0^x(x-s)w(s)u_r(s)\,ds,
\qquad r\ge0.
\end{aligned}
\]
Repeated substitution gives, for $r\ge1$,
\[
\begin{aligned}
u_r(x)=\frac{1}{(a+1)(a+2)}
\int_{0<t_1<\cdots<t_r<x}
t_1^{a+2}(x-t_r)
\times\prod_{i=1}^r w(t_i)
\prod_{i=1}^{r-1}(t_{i+1}-t_i)\,d\mathbf t.
\end{aligned}
\]
In the induction step, the outer integration introduces
$t_{r+1}$ as the largest variable and the new final gap
$x-t_{r+1}$. Integrating the free endpoint
against $(1-x)^b$ gives
\[
\int_{t_r}^1(1-x)^b(x-t_r)\,dx
=\frac{(1-t_r)^{b+2}}{(b+1)(b+2)}.
\]
Comparison with~\eqref{eq:grouped-integral} therefore gives
\[
\int_0^1(1-x)^bu_r(x)\,dx
=\frac{a!b!}{(a+2)!(b+2)!}\,I_r,
\qquad r\ge1.
\]
For $r=0$, the same integral is $a!b!/K!$.

Since $0\le w\le1$, induction gives
$0\le u_r(x)\le x^{2r}/(2r)!$ on $[0,1]$.
Thus the generating function
$u(x,z)=\sum_{r\ge0}u_r(x)z^r$
may be integrated term by term. By~\eqref{eq:gap-bstar},
\begin{equation}\label{eq:S-function}
S(z):=\int_0^1(1-x)^bu(x,z)\,dx
=\frac{2a!b!}{K!}\sum_{r\ge0}
\frac{b_r^*z^r}{(2r+2)!}.
\end{equation}
The recurrence also gives
$u(x,z)=u_0(x)+z\int_0^x(x-s)w(s)u(s,z)\,ds$, or
\begin{equation}\label{eq:volterra}
u_{xx}(x,z)=zw(x)u(x,z)+x^a,
\qquad u(0,z)=u_x(0,z)=0.
\end{equation}

 Next, we construct a symmetric probability measure $\mathcal M$
with Fourier transform $1-t^2S(-t^2)$. The factor $t^2$
will allow us to recover the moments $b_r^*$ by weighting
$\mathcal M$ with $v^2$. The construction uses the fact that
\[
c(x)=\sqrt{w(x)}=x^{a/2}(1-x)^{b/2}
\]
increases up to $p_*=a/(a+b)$ and decreases thereafter.
Lemma~\ref{lem:increasing-speed} can therefore be applied
from either endpoint towards $p_*$.

For real $t$, write $u_t(x)=u(x,-t^2)$ and define
\[
F(x,t)=1-t^2(1-x)^bu_t(x),\qquad 0\le x\le p_*.
\]
Since $w(x)=x^a(1-x)^b$, equation~\eqref{eq:volterra}
now reads $u_t''(x)=x^aF(x,t)$. To see that $F(x,\cdot)$
is a Fourier transform, let
\[
Q(x,t)=\frac{F(x,t)}{(1-x)^b}
=(1-x)^{-b}-t^2u_t(x).
\]
Using $w(x)(1-x)^{-b}=x^a$ in~\eqref{eq:volterra},
we obtain
\begin{equation}\label{eq:gap-Q-equation}
Q_{xx}+t^2wQ=\frac{b(b+1)}{(1-x)^{b+2}},
\qquad Q(0,t)=1,\quad Q_x(0,t)=b.
\end{equation}
Let $h(x,s;t)$ be the homogeneous solution in
Lemma~\ref{lem:increasing-speed}. The function
\[
h(x,0;t)+b\int_0^x\frac{h(x,s;t)}{(1-s)^{b+1}}\,ds
\]
satisfies~\eqref{eq:gap-Q-equation}. Indeed, two
derivatives in $x$, using $h(x,x;t)=1$ and
$h_x(x,x;t)=0$, produce the inhomogeneous term, and the
initial values are $1$ and $b$. Uniqueness gives
\begin{equation}\label{eq:gap-F-mixture}
F(x,t)=(1-x)^b\left[h(x,0;t)
+b\int_0^x\frac{h(x,s;t)}{(1-s)^{b+1}}\,ds\right].
\end{equation}
By the lemma, each $h(x,s;\cdot)$ is the Fourier
transform of a symmetric probability measure. The weights
in~\eqref{eq:gap-F-mixture} are nonnegative and sum to one,
as
\[
(1-x)^b\left(1+b\int_0^x(1-s)^{-b-1}\,ds\right)=1.
\]
Hence $F(x,\cdot)$ is the Fourier transform of a symmetric
probability measure supported on
\[
\left[-\int_0^xc(s)\,ds,\ \int_0^xc(s)\,ds\right].
\]

The corresponding construction from the right endpoint is
obtained by replacing $x$ with $1-x$ and exchanging $a,b$.
Let $v_t$ solve
\[
v_t''+t^2wv_t=(1-x)^b,
\qquad v_t(1)=v_t'(1)=0,
\]
and define
\[
H(x,t)=1-t^2x^av_t(x),\qquad p_*\le x\le1.
\]
Then $v_t''(x)=(1-x)^bH(x,t)$, and $H(x,\cdot)$
is the Fourier transform of a symmetric probability measure
supported on
\[
\left[-\int_x^1c(s)\,ds,\ \int_x^1c(s)\,ds\right].
\]
Both families of measures vary weakly continuously with $x$.

The integral defining $S$ runs over the whole interval,
so the two constructions must be joined at $p_*$. The common
term $-t^2w$ in the equations for $u_t$ and $v_t$ cancels in
\[
\frac{d}{dx}(u_t'v_t-u_tv_t')
=x^av_t-(1-x)^bu_t.
\]
As $v_t(1)=v_t'(1)=0$, integration from $p_*$ to $1$
yields
\begin{equation}\label{eq:gap-parts}
\int_{p_*}^1(1-x)^bu_t(x)\,dx
=\int_{p_*}^1x^av_t(x)\,dx
+u_t'(p_*)v_t(p_*)-u_t(p_*)v_t'(p_*).
\end{equation}
The boundary term can be expressed in terms of $F$ and $H$.
Integrating $u_t''=x^aF$ and $v_t''=(1-x)^bH$
from their respective endpoints gives
\begin{equation}\label{eq:gap-boundary-derivatives}
\begin{aligned}
u_t'(p_*)&=\int_0^{p_*}x^aF(x,t)\,dx,\\
-v_t'(p_*)&=\int_{p_*}^1(1-x)^bH(x,t)\,dx.
\end{aligned}
\end{equation}
Also, by the definitions of $F$ and $H$,
\[
\begin{aligned}
t^2u_t(p_*)&=\frac{1-F(p_*,t)}{(1-p_*)^b},\\
t^2v_t(p_*)&=\frac{1-H(p_*,t)}{p_*^a}.
\end{aligned}
\]
Splitting $S(-t^2)$ at $p_*$, substituting
\eqref{eq:gap-parts} and these four boundary identities,
and multiplying by $t^2$, we obtain
\begin{align}\label{eq:gap-joined}
t^2S(-t^2)
={}&\int_0^{p_*}(1-F(x,t))\,dx
+\int_{p_*}^1(1-H(x,t))\,dx\notag\\
&+(1-H(p_*,t))
\int_0^{p_*}\left(\frac{x}{p_*}\right)^aF(x,t)\,dx\notag\\
&+(1-F(p_*,t))
\int_{p_*}^1\left(\frac{1-x}{1-p_*}\right)^bH(x,t)\,dx.
\end{align}
Since the two intervals partition $[0,1]$, subtracting
\eqref{eq:gap-joined} from one and collecting terms gives
\begin{align}\label{eq:split-positive}
1-t^2S(-t^2)
={}&\int_0^{p_*}F(x,t)\left[
1-\left(\frac{x}{p_*}\right)^a
+\left(\frac{x}{p_*}\right)^aH(p_*,t)
\right]dx\notag\\
&+\int_{p_*}^1H(x,t)\left[
1-\left(\frac{1-x}{1-p_*}\right)^b
+\left(\frac{1-x}{1-p_*}\right)^bF(p_*,t)
\right]dx.
\end{align}

Each bracket in~\eqref{eq:split-positive} is a convex
combination of $1$ and a Fourier transform. Its product
with $F(x,t)$ or $H(x,t)$ is again a Fourier transform,
by convolution of probability measures. Averaging these
measures over $[0,p_*]$ and $[p_*,1]$ gives a symmetric
probability measure $\mathcal M$ with
\[
\widehat{\mathcal M}(t)=1-t^2S(-t^2).
\]
The support of each convolution is also controlled. In the
first integral its radius is at most
\[
\int_0^xc(s)\,ds+\int_{p_*}^1c(s)\,ds
\le\int_0^1c(s)\,ds=L.
\]
In the second integral the corresponding bound is
$\int_x^1c(s)\,ds+\int_0^{p_*}c(s)\,ds\le L$.
Thus $\mathcal M$ is supported on $[-L,L]$.

Comparing coefficients in
$\widehat{\mathcal M}(t)=1-t^2S(-t^2)$
using~\eqref{eq:S-function} gives
\[
\int v^{2r+2}\,\mathcal M(dv)
=\frac{2a!b!}{K!}\,b_r^*,\qquad r\ge0.
\]
The case $r=0$ shows that
\begin{equation}\label{eq:compact-measure}
\mathbb P(Z\in dv)=\frac{K!}{2a!b!}v^2\mathcal M(dv)
\end{equation}
defines a probability law. It is symmetric, supported on
$[-L,L]$, and the same identity gives
$\E Z^{2r}=b_r^*$ for all $r\ge0$.

It remains to identify the factorial factor in
\eqref{eq:gap-factorization}. Since $K=d+1$, grouping the
factors of $(K+dr)!/K!$ by their residues modulo $d$ gives
\begin{align}\label{eq:R-factor}
R_r:=\frac{2(K+dr)!}{K!(a!b!)^r(2r+2)!}
=\left(\frac{d^d}{4a!b!}\right)^r
\frac{\prod_{t=1}^d(1+(t+1)/d)_r}{(3/2)_r(2)_r},
\qquad r\ge0.
\end{align}
Here $(2r+2)!/2=4^r(3/2)_r(2)_r$.
The factor with $t=d-1$ cancels $(2)_r$.
The factor with $t=1$ leaves $(1+2/d)_r/(3/2)_r$,
which is the $r$-th moment of $U$ in~\eqref{eq:gap-Y}
when $d>4$, and equals one when $d=4$. The remaining
factors are the moments of the gamma variables in
\eqref{eq:gap-Y}. Hence $R_r=\E Y^r$. Taking $Y$ independent of $Z$
and using~\eqref{eq:gap-factorization}, we obtain
\[
A_r=R_rb_r^*=\E Y^r\,\E Z^{2r}=\E(YZ^2)^r,
\qquad r\ge0.
\]
\end{proof}

\subsection{The law for $21354$}

We return to the pattern $21354$, used in
Figure~\ref{fig:21354-trunk} to illustrate the
trunk structure. At overlap $m=1$, it has
$a=b=1$ and $d=4$, so the beta factor in
\eqref{eq:gap-Y} is identically one.
Theorem~\ref{thm:gap-two} then reduces to a simpler
representation of its packing counts, as follows.

\begin{corollary}\label{cor:21354}

Let $\nu_s$ be the symmetric probability measures with
Fourier transforms $F(s,t)$ given by~\eqref{eq:gap-F-mixture},
specialized to $a=b=1$, with $0\le s\le1/2$.
Define $\mu$ to be a symmetric probability measure supported on
$[-\pi/8,\pi/8]$ with
\[
\mu(dv)=120v^2\int_0^{1/2}
\bigl[(1-2s)\nu_s
+2s(\nu_{1/2}\ast\nu_s)\bigr](dv)\,ds,
\]
where $\ast$ denotes the convolution of measures.

For independent $Z\sim\mu$, $G_{7/4}$, and
$G_{9/4}$,
\[
N_{21354}^1(r+1)
=\E(64G_{7/4}G_{9/4}Z^2)^r
\qquad(r\ge0).
\]
This Stieltjes moment sequence is determinate.
\end{corollary}

\begin{proof}
When $a=b=1$, we have $p_*=1/2$ and
$H(x,t)=F(1-x,t)$. Thus~\eqref{eq:split-positive},
after substituting $s=1-x$ in its second integral,
is the Fourier transform of
\[
\mathcal M
=2\int_0^{1/2}
\bigl[(1-2s)\nu_s
+2s(\nu_{1/2}*\nu_s)\bigr]\,ds.
\]

The coefficient in~\eqref{eq:compact-measure}
is $5!/2=60$, giving the stated measure.
Its support is contained in $[-L,L]$, where
\[
L=\mathrm B(3/2,3/2)=\frac{\pi}{8}.
\]
Its symmetry and normalization follow from
Theorem~\ref{thm:gap-two}.

Here $d=4$, so $U=1$ in~\eqref{eq:gap-Y},
and the remaining factors give
$Y=64G_{7/4}G_{9/4}$. The moment representation
follows from Theorem~\ref{thm:gap-two}.

Although $Z$ has compact support, this does not
give compact support for the product. We establish
determinacy using the Stieltjes Carleman criterion.

Writing $A_r=N_{21354}^1(r+1)$, the support
bound on $Z$ gives
\[
\begin{aligned}
A_r
&=64^r(7/4)_r(9/4)_r\,\E Z^{2r}\\
&\le \pi^{2r}(7/4)_r(9/4)_r\\
&\le \pi^{2r}(r+2)^{2r},
\qquad r\ge1.
\end{aligned}
\]
The last inequality holds because every factor
in $(7/4)_r$ and $(9/4)_r$ is at most $r+2$.
Therefore
\[
\sum_{r\ge1}A_r^{-1/(2r)}
\ge\frac1\pi\sum_{r\ge1}\frac1{r+2}
=\infty.
\]
By the Stieltjes Carleman criterion
\cite[Theorem~1.11]{ST43},
the sequence is determinate.
\end{proof}

\section{Explicit laws and a Jacobi continued fraction}\label{sec:explicit}

The preceding sections establish moment representations for several
families, often through explicit product or integral constructions. We conclude with two examples that require different methods. For $135246$ at overlap three, a recurrence for shifted tableaux
leads to a compactly supported representing measure. While Theorem~\ref{thm:grid} establishes Stieltjes positivity for two-level
shifted grids without generally identifying their representing
measures, here we determine the
Jacobi continued fraction for $14325$ at overlap two using weighted paths.

\subsection{A compactly supported law for 135246 at $m=3$}

\begin{theorem}\label{thm:135246}
Let $A_r=N_{135246}^3(r+1)$. Let $U,V,W$ be independent, with
\[
U\sim\tfrac9{16}\operatorname{Beta}(1/2,2)
       +\tfrac7{16}\operatorname{Beta}(3/2,1),
\]
\[
V\sim\operatorname{Beta}(4/3,5/3),\qquad
W\sim\operatorname{Beta}(5/3,4/3),
\]
Write $\nu=\Law(27\,UVW)$.
Then
\begin{equation}\label{eq:135-law}
\mu(dx)=\frac{2x}{1+x}\,\nu(dx),\quad\quad x\in [0,27],
\end{equation}
is the unique  probability measure on $\mathbb R$ with moments $A_r$.
\end{theorem}

\begin{proof}

\medskip\noindent\emph{(1) From the grid to a recurrence.}

For $135246$, the prefix $(1,3,5)$ and suffix $(2,4,6)$
give a three-level shifted grid in which all subdivision
parameters vanish (Section~\ref{sec:grids}).
In turn, as shown in Figure~\ref{fig:135246-shifted}, the map
$u_{i,q}\mapsto(i,q+i)$ identifies this grid with a
`shifted strip' in the sense of \cite{Sun2017}.
Writing $g_{3,n}$ for the number of standard tableaux of this
three-row strip with $n$ cells per row, it immediately follows that:
\[
N_{135246}^3(j)=g_{3,j+1}\quad(j\ge1),
\qquad A_r=g_{3,r+2}.
\]

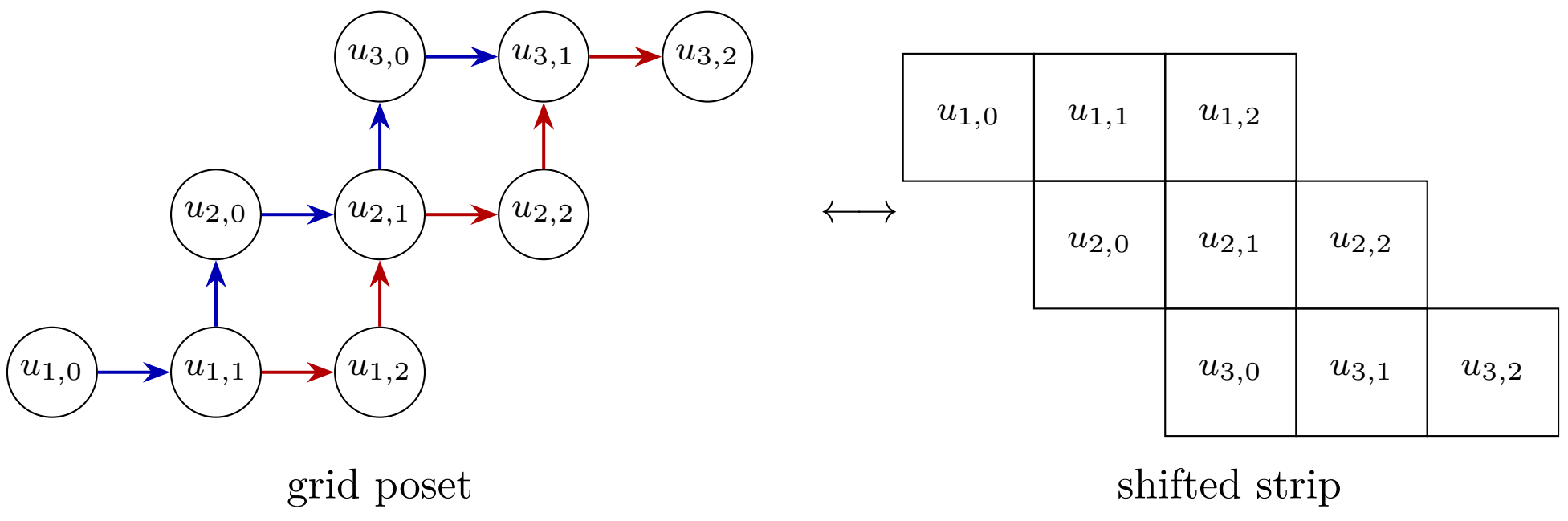
\begin{figure}[htbp]
\centering
\begin{tikzpicture}[>=Stealth,font=\small,
  vertex/.style={circle,draw,fill=white,inner sep=1.9pt,font=\scriptsize},
  first/.style={->,draw=blue!70!black,line width=.75pt},
  second/.style={->,draw=red!70!black,line width=.75pt},
  tabcell/.style={draw,minimum width=10.8mm,minimum height=10.5mm,
                 inner sep=0pt,font=\scriptsize}]

\foreach \i in {1,2,3}{
  \foreach \q in {0,1,2}{
    \pgfmathsetmacro{\xx}{1.35*(\q+\i-1)}
    \pgfmathsetmacro{\yy}{1.30*(\i-1)}
    \node[vertex] (g\i\q) at (\xx,\yy) {$u_{\i,\q}$};
  }
}

\foreach \i in {1,2,3}{
  \draw[first]  (g\i0) -- (g\i1);
  \draw[second] (g\i1) -- (g\i2);
}
\foreach \i/\ii in {1/2,2/3}{
  \draw[first]  (g\i1) -- (g\ii0);
  \draw[second] (g\i2) -- (g\ii1);
}

\foreach \i in {1,2,3}{
  \foreach \q in {0,1,2}{
    \pgfmathsetmacro{\xx}{7.55+1.08*(\q+\i-1)}
    \pgfmathsetmacro{\yy}{2.10-1.05*(\i-1)}
    \node[tabcell] at (\xx,\yy) {$u_{\i,\q}$};
  }
}

\node at (6.65,1.30) {$\longleftrightarrow$};
\node[font=\small] at (2.70,-.95) {grid poset};
\node[font=\small] at (9.70,-.95) {shifted strip};
\end{tikzpicture}
\caption{The two-occurrence grid for $135246$ and its `shifted strip',
under $u_{i,q}\mapsto(i,q+i)$. Edge colors identify distinct occurrences.}
\label{fig:135246-shifted}
\end{figure}

The tableau interpretation allows us to use an enumeration result of Sun \cite{Sun2017}, who studied shifted strips with equal row lengths and derived recurrences for their tableau counts. In the three-row case relevant here, the recurrence takes the form:
\[
g_{3,n+1}+g_{3,n}
=\frac{7n+1}{n^2(n+1)^2}
  \binom{2n-2}{n-1}\binom{3n}{n-1}
\qquad(n\ge1).
\]
Here, $g_{3,1}=1$ is the auxiliary initial value (and not the extension count for three disconnected single-cell rows). The one-occurrence grid is a chain, so $g_{3,2}=1$ as expected.

Since $A_r=g_{3,r+2}$, Sun's recurrence becomes
\begin{equation}\label{eq:Sun-recurrence}
A_r+A_{r-1}=\tau_{r+1}\qquad(r\ge1),
\end{equation}
where
\begin{equation}\label{eq:135246-tau}
\tau_n=\frac{(7n+1)(3n)!}
{2(2n-1)(2n+1)(n+1)^2(n!)^3}
\qquad(n\ge1),
\end{equation}
with $A_0=1$.

\medskip\noindent\emph{(2) From the recurrence to moments.}

The mixture defining $U$ has moments
\[
\mathbb E U^r
=\frac{7r+8}{8}\frac{(1/2)_r}{(5/2)_r},
\]
where $(\cdot)_r$ once again denotes the rising factorial, i.e.
$(a)_r:=a(a+1)\cdots(a+r-1)$, with $(a)_0=1$.

Let $X=27\,UVW$. Recalling that $U,V,W$ are independent and writing out the beta moments for $V$ and $W$, we obtain
\begin{equation}\label{eq:135-beta-moments}
2\mathbb E X^r
=2\cdot27^r\frac{7r+8}{8}
 \frac{(1/2)_r(4/3)_r(5/3)_r}
      {(5/2)_r(3)_r^2}
=\tau_{r+1}\qquad(r\ge0),
\end{equation}
where the last equality follows by simplifying the
factorials. 

Consider the positive finite measure
\[
\mu(dx)=\frac{2x}{1+x}\,\nu(dx),
\qquad \nu=\operatorname{Law}(X).
\]
Its $r$th moment is $2\,\mathbb E[X^{r+1}/(1+X)]$. The sum of two consecutive moments of $\mu$ is therefore
\[
2\,\mathbb E\left[
\frac{X^{r+1}+X^r}{1+X}
\right]
=2\,\mathbb E X^r=\tau_{r+1} \qquad(r\ge1).
\]
Observe that this is the same recurrence as \eqref{eq:Sun-recurrence}.

Since the moments of $\mu$ satisfy the same recurrence as $A_r$,
all that remains is to establish agreement at $r=0$, namely
\[
A_0=2\,\mathbb E\frac{X}{1+X}.
\]
Since $A_0=1$ and $X/(1+X)=1-1/(1+X)$,
we are equivalently looking to show that
\begin{equation}
2\,\mathbb E\frac{1}{1+X}=1.\label{eq-exp-1}
\end{equation}

\medskip\noindent\emph{(3) The normalization.}

To prove \eqref{eq-exp-1}, we start by setting 
\[
F(t)=1-2t\,\mathbb E\frac{1}{1+tX},\quad \quad t\ge0.
\]
To evaluate $F(1)$, we first identify the Taylor coefficients
of $F$ at the origin. By~\eqref{eq:135-beta-moments},
\[
F(0)=1,\qquad [t^n]F(t)=(-1)^n\tau_n
\quad(n\ge1).
\]

We use an integral ansatz whose Taylor coefficients can be
checked against these values and whose integral at $t=1$
can be evaluated explicitly. For $t\ge0$ and $0\le u\le1$,
let $y_t(u)\in(0,1]$ be the unique solution of
\begin{equation}\label{eq:135246-cubic}
y_t(u)+4tu\,y_t(u)^3=1.
\end{equation}
(Existence and uniqueness follow since $y\mapsto y+4tu\,y^3$
is strictly increasing on $[0,1]$, with values $0$ and
$1+4tu\ge1$ at the endpoints.)
We claim that
\begin{equation}\label{eq:135246-integral}
F(t)=\frac1\pi\int_0^1
\sqrt{\frac{1-u}{u}}\,
\bigl(1+2y_t(u)-y_t(u)^2-4t\,y_t(u)^3\bigr)\,du.
\end{equation}

We first verify this identity near $t=0$.
Lagrange inversion applied to~\eqref{eq:135246-cubic}
gives, for $n\ge1$,
\[
\begin{aligned}
&[t^n]\bigl(1+2y_t(u)-y_t(u)^2-4t\,y_t(u)^3\bigr)=
(-4)^n\frac{(3n)!}{n!(2n+1)!}
\left(\frac{1-n}{n+1}u^n+u^{n-1}\right).
\end{aligned}
\]
Together with the beta integral
\[
\frac1\pi\int_0^1
u^{m-\frac12}(1-u)^{\frac12}\,du
=\frac{1}{2\cdot4^m(m+1)}\binom{2m}{m}
\qquad(m\ge0),
\]
we find that the coefficient of $t^n$ on the right-hand
side of~\eqref{eq:135246-integral} is
\[
\begin{aligned}
\frac{(-1)^n(3n)!}{2(2n+1)(n!)^3}
\left(\frac{1-n}{(n+1)^2}+\frac{2}{2n-1}\right)
=
\frac{(-1)^n(7n+1)(3n)!}
{2(2n-1)(2n+1)(n+1)^2(n!)^3}
=(-1)^n\tau_n,
\end{aligned}
\]
where the final equality uses~\eqref{eq:135246-tau}.

At $t=0$, the right-hand side of~\eqref{eq:135246-integral}
equals $1=F(0)$. Together with the agreement of the
coefficients of $t^n$ for $n\ge1$, this establishes
\eqref{eq:135246-integral} in a neighborhood of $t=0$.

Since $0\le X\le27$, $F$ is real analytic for $t\ge0$.
The integral on the right-hand side
of~\eqref{eq:135246-integral} is real analytic as well.
Indeed, the derivative with respect to $y$ in
\eqref{eq:135246-cubic} is
$1+12tu\,y_t(u)^2\ge1$.
The analytic implicit-function theorem, together with
compactness of $[0,1]$, gives analytic dependence
of $y_t(u)$ on $t$, locally uniformly in $u$.
Since the weight $\sqrt{(1-u)/u}$ is integrable,
this analyticity is preserved under integration.

Both sides of~\eqref{eq:135246-integral} are therefore
real analytic for $t\ge0$. Since they agree near $t=0$,
the identity theorem establishes
\eqref{eq:135246-integral} for all $t\ge0$,
in particular at $t=1$.
It remains to evaluate the integral
in~\eqref{eq:135246-integral} at $t=1$.

Write $y(u)=y_1(u)$ and $h(u)=\sqrt{u(1-u)}$.
We claim that the function under the integral sign
has the antiderivative
\begin{equation}\label{eq:135246-primitive}
G(u)=(4y^3-4y^2+3y+1)h
      -4\arcsin(2y^2h).
\end{equation}
To verify \eqref{eq:135246-primitive}, we differentiate the two terms.
At $t=1$, equation~\eqref{eq:135246-cubic} becomes
$y+4uy^3=1$. Differentiating with respect to $u$
and using this equation once more gives
\[
(1+12uy^2)y'=-4y^3,
\qquad
y'=-\frac{4y^4}{3-2y}.
\]
For the arcsine term, the chain rule introduces
the square root $\sqrt{1-4y^4h^2}$.
The cubic also gives $u=(1-y)/(4y^3)$, so
$h^2=u(1-u)$ implies
\[
\sqrt{1-4y^4h^2}
=\frac{2y^2-y+1}{2y}>0.
\]
Together with $h'=(1-2u)/(2h)$, this gives
\[
\begin{aligned}
\frac{d}{du}\arcsin(2y^2h)
=\frac{4yy'h+2y^2h'}{\sqrt{1-4y^4h^2}}
=\frac{2y^2-1}{(3-2y)h}.
\end{aligned}
\]

We can now differentiate~\eqref{eq:135246-primitive}.
The product rule and the identities above give
\[
\begin{aligned}
G'(u)
&=(12y^2-8y+3)y'h
  +(4y^3-4y^2+3y+1)h'-\frac{4(2y^2-1)}{(3-2y)h}\\
&=\sqrt{\frac{1-u}{u}}\,
  (1+2y-y^2-4y^3),
\end{aligned}
\]
where the last equality follows by substituting
the expressions for $y'$ and $h'$ and simplifying
with $4uy^3=1-y$ and $h^2=u(1-u)$.
This is precisely the function under the integral
sign in~\eqref{eq:135246-integral} at $t=1$,
as claimed.

Finally, $h(0)=h(1)=0$ and $y$ is bounded.
Both terms in~\eqref{eq:135246-primitive} therefore
vanish at the endpoints, the second because
$\arcsin(0)=0$. Hence $G(0)=G(1)=0$.
Since the integrand is integrable at both endpoints,
we conclude that
\[
F(1)=\frac1\pi\int_0^1G'(u)\,du
=\frac{G(1)-G(0)}{\pi}=0,
\]
which proves~\eqref{eq-exp-1}.  Since $\mu$ is compactly 
supported, its moment sequence is determinate.
\end{proof}

\subsection{The Jacobi fraction for $14325$ at $m=2$}

We establish the following.

\begin{theorem}\label{thm:14325}
Let $A_r=N_{14325}^2(r+1)$ for $r\ge0$. Then, as formal power series,
\begin{equation}\label{eq:14325-fraction}
\sum_{r\ge0}A_rz^r
=\cfrac1{1-6z-\cfrac{35z^2}{1-18z-\cfrac{143z^2}{1-30z-\cfrac{323z^2}{1-42z-\ddots}}}},
\end{equation}
that is, with Jacobi coefficients 
\[
\alpha_n=12n+6\quad(n\ge0)\qquad\text{and}\qquad
\beta_n=36n^2-1\quad(n\ge1).
\]
\end{theorem}

The proof relies on a special structural feature of
$14325$ at overlap two. Concretely, inserting another occurrence into a packing gives a
recurrence for the packing counts, which we interpret in terms
of weighted paths with steps $+1$ and $-2$. Since both steps
increase the height by one modulo three, grouping them into
triples produces a four-diagonal (infinite) matrix $Q$ satisfying
$$A_r=(Q^r)_{0,0}=e_0^{T}Q^r e_0$$
where $e_0=(1,0,\ldots)^{T}$. This is not yet a Jacobi representation, as the contracted paths can still descend by two levels. However, we exhibit an explicit unit lower-triangular
matrix $R$ that algebraically conjugates $Q$ to a tridiagonal Jacobi matrix $J$.
Moreover, $R$ satisfies $Re_0=e_0$ and
$e_0^{T}R=e_0^{T}$, preserving the moments.

\begin{lemma}\label{lem:14325-paths}
For $j\ge1$, the count $N_{14325}^2(j)$ equals the total weight of
nonnegative paths in $\mathbb Z_{\geq 0}\times \mathbb Z_{\geq 0}$ from height zero to height zero with:\begin{itemize}\item $2j$ up-steps,
$U=(1,1)$, of weight one, and 
\item $j$ down-steps, $D=(1,-2)$, with a down-step starting at height $h$ carrying weight $h-1$.
\end{itemize}
\end{lemma}

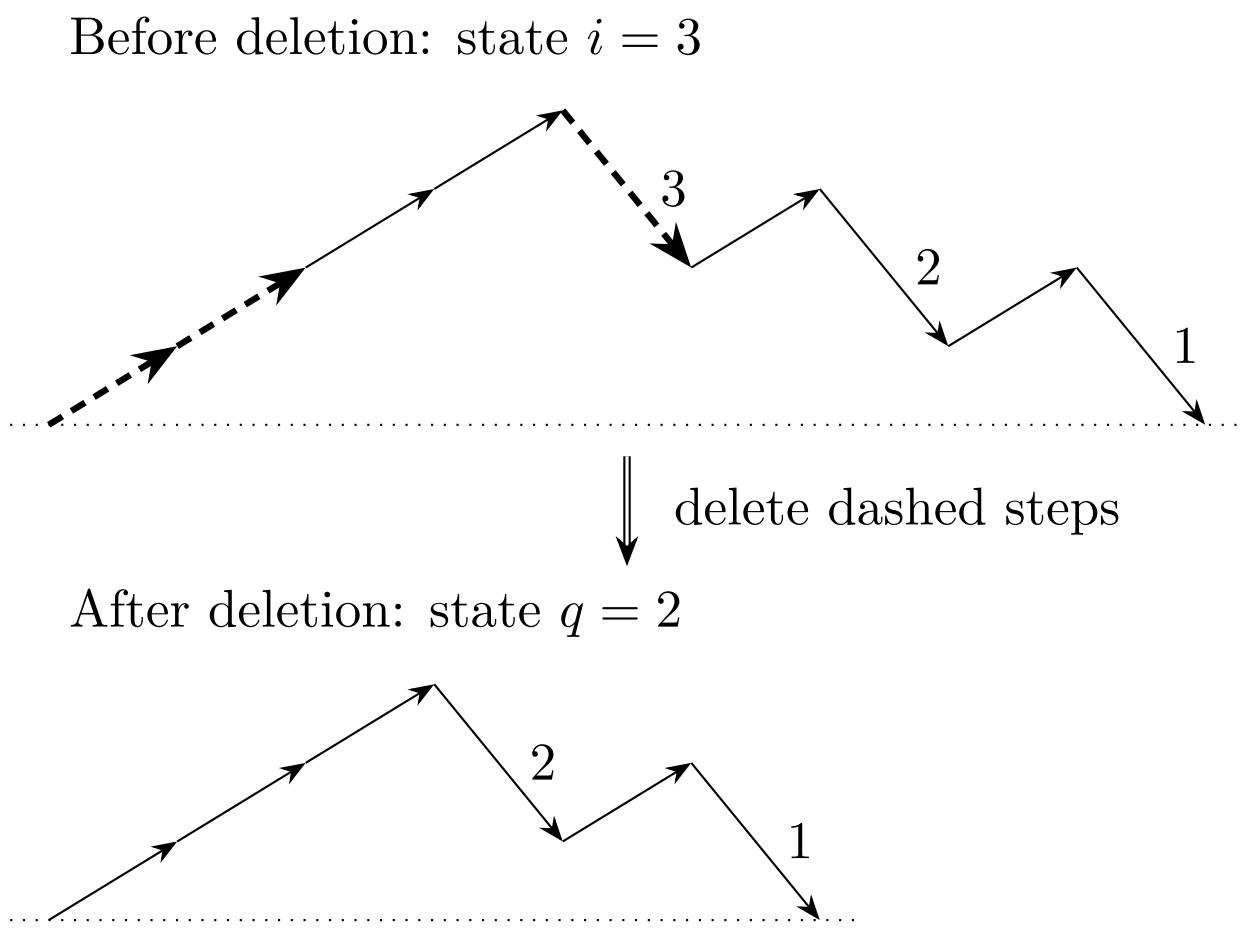
\begin{figure}[t]
 \centering
 \begin{tikzpicture}[>=Stealth,font=\small,x=8.5mm,y=5.2mm]
  \node[anchor=west] at (0,11.25) {Before deletion: state $i=3$};
  \draw[dotted] (-.3,6.3)--(9.3,6.3);
  \draw[->,densely dashed,very thick] (0,6.3)--(1,7.3);
  \draw[->,densely dashed,very thick] (1,7.3)--(2,8.3);
  \draw[->] (2,8.3)--(3,9.3);
  \draw[->] (3,9.3)--(4,10.3);
  \draw[->,densely dashed,very thick] (4,10.3)--node[right=2pt] {$3$}(5,8.3);
  \draw[->] (5,8.3)--(6,9.3);
  \draw[->] (6,9.3)--node[right=2pt] {$2$}(7,7.3);
  \draw[->] (7,7.3)--(8,8.3);
  \draw[->] (8,8.3)--node[right=2pt] {$1$}(9,6.3);

  \draw[->,double] (4.5,5.9)--(4.5,4.5);
  \node[right=5pt] at (4.5,5.2) {delete dashed steps};

  \node[anchor=west] at (0,3.9) {After deletion: state $q=2$};
  \draw[dotted] (-.3,0)--(6.3,0);
  \draw[->] (0,0)--(1,1);
  \draw[->] (1,1)--(2,2);
  \draw[->] (2,2)--(3,3);
  \draw[->] (3,3)--node[right=2pt] {$2$}(4,1);
  \draw[->] (4,1)--(5,2);
  \draw[->] (5,2)--node[right=2pt] {$1$}(6,0);
\end{tikzpicture}
  \caption{Deleting the two initial up-steps and the first down-step (dashed) from a path of state $i=3$ leaves a path of state $q=2$. The weights $2$ and $1$ of the remaining down-steps are unchanged. Conversely, insertion multiplies the path weight by $3$.} \label{fig:14325-path-insertion}
  \end{figure}

\begin{proof}
We show that the linear extensions of the packing
poset satisfy the same recurrence as these weighted paths.

By the shifted-grid construction of Section~\ref{sec:grids},
the overlap-two packing poset for $14325$ has two increasing
chains $x_0<\cdots<x_j$ and $y_0<\cdots<y_j$,
together with vertices $z_0,\ldots,z_{j-1}$. For $0\le q<j$, the $(q+1)$-st
occurrence consists, in positional order, of
$x_q,y_q,z_q,x_{q+1},y_{q+1}$.
Since the pattern is $14325$, the additional relations are
\begin{equation}\label{eq:14325-orders}
x_{q+1}<z_q<y_q
\qquad(0\le q<j).
\end{equation}
Together with the two previously specified chains, these relations
generate the entire packing poset. 

Call an extension of this poset in \emph{state $i$} if exactly $i+2$
vertices follow $x_j$, and let $b_j(i)$ count such extensions.
Since $x_j<z_{j-1}<y_{j-1}<y_j$, at least three
vertices follow $x_j$ in every linear extension,
so $i\ge1$.
For $j=1$, the poset is a chain, so $b_1(1)=1$ and
$b_1(i)=0$ for $i>1$.

To pass from $j$ to $j+1$ occurrences, insert
$x_{j+1},z_j,y_{j+1}$ subject to
\[
x_j<x_{j+1}<z_j<y_j<y_{j+1}.
\]
The old last vertex is $y_j$, and the new last vertex must be
$y_{j+1}$. Suppose the old extension has state $q$ and the new
one has state $i$. Since $z_j$ and $y_{j+1}$ lie after
$x_{j+1}$, precisely $i$ old vertices must remain after
$x_{j+1}$. Its insertion position is therefore fixed, and is
available exactly when $1\le i\le q+2$. Of these $i$ old
vertices, the last is $y_j$, leaving $i$ possible positions
for $z_j$ before $y_j$. Deleting the three new vertices reverses
the construction. Consequently,
\begin{equation}\label{eq:14325-insertion}
b_{j+1}(i)=i\sum_{q\ge i-2}b_j(q),
\qquad i\ge1,
\end{equation}
where $b_j(q)=0$ for $q<1$.

We now show that the same recurrence governs the weighted paths. For an illustration, see Figure~\ref{fig:14325-path-insertion}.

A nonnegative
path has at least two up-steps before its first down-step. We
assign it state $i$ when that down-step starts at height $i+1$.
The only path with one down-step is $UUD$, of state one and
weight one.

Given a path with $j+1$ down-steps, delete its first two
up-steps and its first down-step. The resulting path is
nonnegative and has $j$ down-steps. The deleted steps have
net height change zero, so the starting heights, and hence
the weights, of the remaining down-steps are unchanged.

Conversely, suppose a path with $j$ down-steps has state $q$,
so it begins with $q+1$ up-steps. Prepend two up-steps and
insert a down-step after the first $i-1$ original up-steps,
where $1\le i\le q+2$. The inserted down-step starts at
height $i+1$, has weight $i$, and becomes the first down-step.
These operations are inverse. Thus the total path weights
in state $i$ satisfy~\eqref{eq:14325-insertion} with the same
initial condition as $b_j(i)$. Summing over $i$ proves the lemma.
\end{proof}

\begin{lemma}\label{lem:14325-contraction}
Let $Q$ be the matrix indexed by nonnegative integers, with
nonzero entries
\begin{equation}\label{eq:Q14325}
\begin{aligned}
Q_{n,n+1}&=1,&Q_{n,n}&=9n+6,\\
Q_{n,n-1}&=27n^2+9n-1 &&(n\ge1),\\
Q_{n,n-2}&=(3n-3)(3n-1)(3n+1) &&(n\ge2).
\end{aligned}
\end{equation}
Then $A_r=(Q^r)_{0,0}$ for all $r\ge0$.
\end{lemma}

\begin{proof}
By Lemma~\ref{lem:14325-paths}, $A_r$ is the total weight
of paths with $r+1$ down-steps. Every such path begins with
two up-steps and ends with a down-step from height two.
These three steps have weight one. Removing them leaves
a nonnegative path of $3r$ steps from height two to height two.

Both $U$ and $D$ increase the height by one modulo three.
Therefore the height at each boundary between successive
triples is $3n+2$ for some $n\ge0$. We take $n$ as the
contracted state. The possible triples and their weights,
starting at height $3n+2$, are
\[
\begin{array}{c|c|l}
\text{state change}&\text{triple}&\text{weight}\\\hline
+1&UUU&1\\
0&DUU&3n+1\\
 &UDU&3n+2\\
 &UUD&3n+3\\
-1&DDU&(3n+1)(3n-1)\\
 &DUD&(3n+1)(3n)\\
 &UDD&(3n+2)(3n)\\
-2&DDD&(3n+1)(3n-1)(3n-3).
\end{array}
\]
The triples decreasing the state by one require $n\ge1$, while 
the triple $DDD$ requires $n\ge2$. (At lower heights these triples would take the path below zero.) All other triples in the table stay
nonnegative.

Summing the weights for each possible state change gives
the entries of $Q$ in~\eqref{eq:Q14325}. Conversely,
each contracted transition can be expanded uniquely once
its constituent triple is chosen, with the same weight.
The paths counted by $A_r$ therefore correspond to
length-$r$ weighted walks on the nonnegative states,
starting and ending at zero. Their total weight is
$(Q^r)_{0,0}$. For $r=0$, both sides equal one.
\end{proof}

\begin{figure}[htbp]
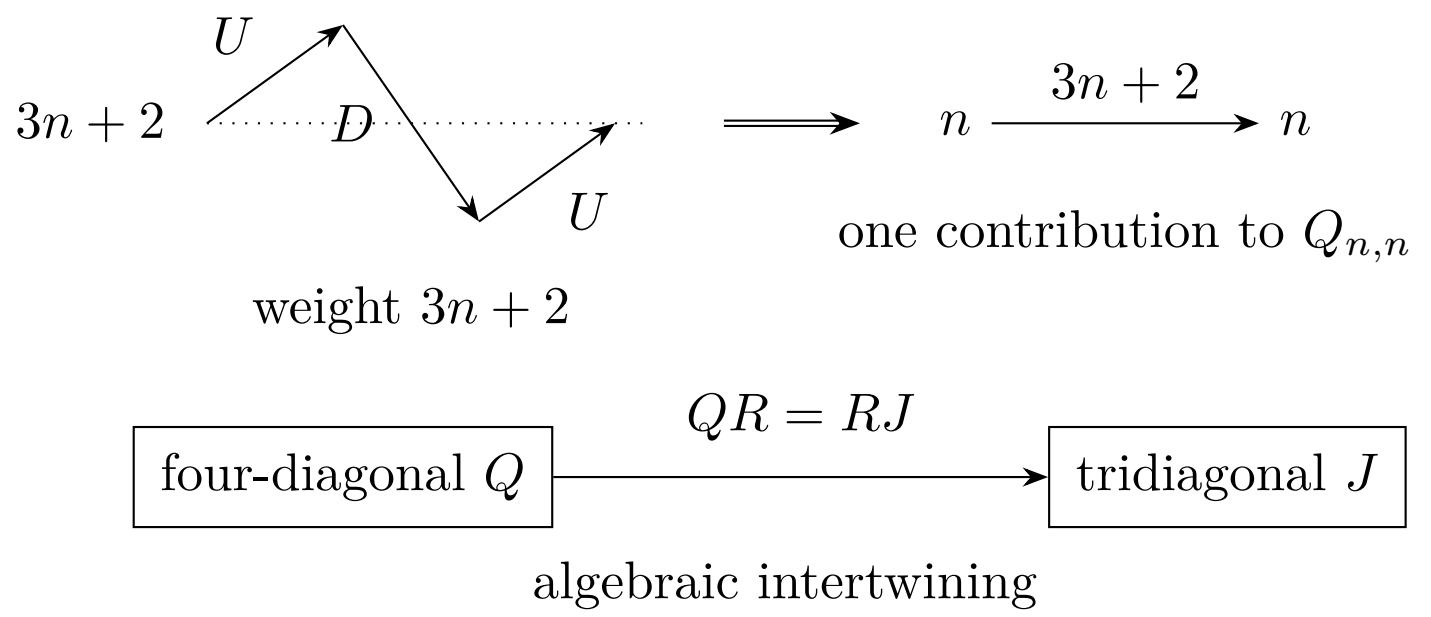

\centering\ContractionFigure
\caption{A triple $UDU$ starts and ends at height $3n+2$.
Its down-step has weight $3n+2$, one of the three contributions
to $Q_{n,n}$. The triple $DDD$ accounts for the second
subdiagonal of $Q$.}
\label{fig:contraction}
\end{figure}

The nontrivial second subdiagonal, $Q_{n,n-2}$, prevents $Q$ from being a Jacobi matrix.
We now exhibit an explicit intertwining $QR=RJ$ with a
tridiagonal Jacobi matrix $J$, and show that $Q$ and $J$
have the same $(0,0)$ moments.

\begin{lemma}\label{lem:14325-intertwining}
Define $R$ by
\begin{equation}\label{eq:R14325}
R_{n,k}=\frac{n!}{(2k-n)!(n-k)!}
       \left(\frac32\right)^{n-k}
\quad\text{if }\left\lceil\frac n2\right\rceil\le k\le n,
\end{equation}
and $R_{n,k}=0$ otherwise. Let $J$ be the tridiagonal matrix
with nonzero entries
\begin{equation}\label{eq:J14325}
J_{n,n+1}=1,\quad J_{n,n}=12n+6\quad(n\ge0),
\qquad J_{n,n-1}=36n^2-1\quad(n\ge1).
\end{equation}
Then $QR=RJ$.
\end{lemma}

\begin{proof}
We verify the matrix identity entrywise. All sums are finite,
since $Q$ and $J$ are banded and each row of $R$ has finite
support. With entries involving negative indices taken to be
zero, the two sides are
\begin{align}
(QR)_{n,k}={}&R_{n+1,k}+(9n+6)R_{n,k}
 +(27n^2+9n-1)R_{n-1,k}\notag\\
&+(3n-3)(3n-1)(3n+1)R_{n-2,k},
\label{eq:QR-entry}\\
(RJ)_{n,k}={}&R_{n,k-1}+(12k+6)R_{n,k}
 +[36(k+1)^2-1]R_{n,k+1}.
\label{eq:RJ-entry}
\end{align}

First suppose $n\ge2$ and
$\lceil n/2\rceil\le k\le n$. Set
\[
p=n-k,\qquad s=2k-n,\qquad c=R_{n,k}>0,
\]
so that $n=s+2p$ and $k=s+p$. From~\eqref{eq:R14325},
\[
\begin{aligned}
\frac{R_{n+1,k}}c&=\frac{3s(n+1)}{2(p+1)},&
\frac{R_{n,k-1}}c&=\frac{3s(s-1)}{2(p+1)},\\
\frac{R_{n-1,k}}c&=\frac{2p}{3n(s+1)},&
\frac{R_{n-2,k}}c&=\frac{4p(p-1)}{9n(n-1)(s+1)(s+2)},\\
\frac{R_{n,k+1}}c&=\frac{2p}{3(s+1)(s+2)}.&&
\end{aligned}
\]
When a numerator is zero, the corresponding entry of $R$
lies outside its support, so these ratios still apply.

Substitution into~\eqref{eq:QR-entry} and~\eqref{eq:RJ-entry}
gives
\[
\begin{aligned}
\frac{(QR-RJ)_{n,k}}c
&=\frac{3s[(n+1)-(s-1)]}{2(p+1)}+9n-12k\\
&\quad+\frac{2p}{3}\left[
\frac{27n^2+9n-1}{n(s+1)}
+\frac{2(p-1)(9n^2-1)}{n(s+1)(s+2)}
-\frac{36(k+1)^2-1}{(s+1)(s+2)}\right].
\end{aligned}
\]
Using $n=s+2p$ and $k=s+p$, the first line is $6p$
and the bracket equals $-9$. For the latter, multiplying
by $n(s+1)(s+2)$ reduces the claim to
\[
\begin{aligned}
0={}&(27n^2+9n-1)(s+2)+2(p-1)(9n^2-1)\\
&-n[36(k+1)^2-1]+9n(s+1)(s+2),
\end{aligned}
\]
which is the identity obtained upon substituting
$s=n-2p$ and $k=n-p$. Thus
$(QR-RJ)_{n,k}/c=6p-6p=0$ throughout the support of $R$.

It remains to check the entries adjacent to that support.
For $k=n+1$, both~\eqref{eq:QR-entry} and
\eqref{eq:RJ-entry} equal one. At the lower boundary,
the possible extra entry depends on the parity of $n$.

If $n=2p\ge2$ and $k=p-1$, only the last term in
\eqref{eq:QR-entry} and the last term in
\eqref{eq:RJ-entry} can be nonzero. The identities
\[
R_{2p,p}=3(2p-1)R_{2p-2,p-1},\qquad
(3n-3)(3n-1)(3n+1)=3(2p-1)(36p^2-1)
\]
show that both sides equal $(36p^2-1)R_{2p,p}$.

If $n=2p+1\ge3$ and $k=p$, let
$B=R_{2p+1,p+1}$. Then
\[
R_{2p,p}=\frac{B}{2p+1},\qquad
R_{2p-1,p}=\frac{B}{3(2p+1)}.
\]
The entries in~\eqref{eq:QR-entry} and
\eqref{eq:RJ-entry} therefore satisfy
\[
\begin{aligned}
\frac{(QR)_{2p+1,p}}B
&=\frac{27n^2+9n-1+\frac13(3n-3)(3n-1)(3n+1)}{n}\\
&=9(n+1)^2-1=36(p+1)^2-1
=\frac{(RJ)_{2p+1,p}}B.
\end{aligned}
\]
The first two rows give, directly,
\[
(QR)_{0,k}=(RJ)_{0,k}=(6,1,0,\ldots)_k,
\qquad
(QR)_{1,k}=(RJ)_{1,k}=(35,18,1,0,\ldots)_k.
\]
All entries further outside the support vanish on both
sides. This proves $QR=RJ$.
\end{proof}

\begin{proof}[Proof of Theorem~\ref{thm:14325}]
By Lemmas~\ref{lem:14325-contraction} and
\ref{lem:14325-intertwining},
$A_r=(Q^r)_{0,0}$ and $QR=RJ$. Since $Q$ and $J$
are banded and the rows and columns of $R$ have finite
support, all matrix products below involve finite sums
at each fixed power. Induction gives
\[
Q^rR=RJ^r\qquad(r\ge0).
\]
Formula~\eqref{eq:R14325} also shows that
$Re_0=e_0$ and $e_0^{T}R=e_0^{T}$,
where $e_0$ is the first coordinate vector. Hence
\[
A_r=e_0^{T}Q^rRe_0
=e_0^{T}RJ^re_0
=(J^r)_{0,0}.
\]

We can now use the weighted-path interpretation of a
tridiagonal matrix. By~\eqref{eq:J14325}, $(J^r)_{0,0}$
counts weighted Motzkin excursions of length $r$, where an up-step
has weight one, a level-step at height $n$ has weight
$\alpha_n=12n+6$, and a down-step from height $n$
has weight $\beta_n=36n^2-1$. By Flajolet's correspondence between weighted Motzkin paths
and Jacobi continued fractions~\cite[Theorem~1]{Flajolet1980},
the generating function of these excursions has the claimed
Jacobi fraction.
\end{proof}

We conclude by extending the path interpretation of
Lemma~\ref{lem:14325-paths} to a larger family of patterns.
The resulting generating functions admit branched Stieltjes
continued fractions of the type considered by Albenque and
Bouttier~\cite{AlbenqueBouttier2012} and subsequently
developed by P\'etr\'eolle, Sokal and
Zhu~\cite{PetreolleSokalZhu2023}.

\begin{proposition}\label{prop:14325-family}
Let $b,c\ge0$, $K=b+c+4$ and $e=b+c+1$, and
suppose that $\pi\in S_K$ satisfies
\[
(\pi_1,\pi_2)=(1,b+3),
\qquad
(\pi_{K-1},\pi_K)=(2,K).
\]
For $h\ge0$, let $G_h(z)$ be the formal power series
with constant term one determined by
\begin{equation}\label{eq:branched-family}
G_h(z)=
\frac{1}{
1-\binom{h+b}{b}z
\prod_{\ell=1}^{e}G_{h+\ell}(z)}.
\end{equation}
Then
\begin{equation}\label{eq:branched-packing}
\sum_{r\ge0}N_\pi^2(r+1)z^r
=\prod_{\ell=0}^{e}G_\ell(z).
\end{equation}
\end{proposition}

\begin{proof}
The packing poset has two increasing chains
$x_0<\cdots<x_j$ and $y_0<\cdots<y_j$. Each occurrence
adds a chain of $b$ vertices between $x_{q+1}$ and
$y_q$, and a chain of $c$ vertices between $y_q$ and
$y_{q+1}$.

As in the proof of Lemma~\ref{lem:14325-paths},
we count extensions according to the number of vertices
following $x_j$. Let $a_j(i)$ count those in which
exactly $i+e$ vertices follow $x_j$. Thus
$a_1(1)=1$, and $a_1(i)=0$ otherwise.

Suppose an extension with $j$ occurrences has state
$k$. When adding another occurrence, the new vertex
$x_{j+1}$ can be inserted so that exactly $i$ old
vertices follow it, for any $1\le i\le k+e$.
The last of these vertices is $y_j$. The new chain
of $b$ vertices can be interleaved with the other
$i-1$ old vertices in
$\binom{i+b-1}{b}$ ways. The remaining new vertices
are forced after $y_j$. Consequently,
\begin{equation}\label{eq:branched-insertion}
a_{j+1}(i)=
\binom{i+b-1}{b}
\sum_{k\ge i-e}a_j(k),
\qquad i\ge1.
\end{equation}

Consider now nonnegative paths with up-steps $+1$,
of weight one, and down-steps $-e$, a down-step
from height $h$ having weight
\[
w(h)=\binom{h-c-1}{b}.
\]
Assign a path state $i$ when its first down-step
starts at height $i+e-1$. That step has weight
$\binom{i+b-1}{b}$. Deleting the first $e$ up-steps
and the first down-step gives the same recurrence
\eqref{eq:branched-insertion}, exactly as in
Lemma~\ref{lem:14325-paths}. The initial path
$U^eD$ has weight one. Hence $N_\pi^2(j)$ equals
the total weight of excursions with $j$ down-steps.

Let $G_h(z)$ count excursions starting and ending
at height $h$ without visiting a lower height,
with $z$ marking down-steps. The first-return
decomposition gives
\[
G_h(z)=1+w(h+e)z
\prod_{\ell=0}^{e}G_{h+\ell}(z).
\]
Since $w(h+e)=\binom{h+b}{b}$, this is equivalent
to \eqref{eq:branched-family}. Finally,
\[
G_0(z)=1+\sum_{j\ge1}N_\pi^2(j)z^j.
\]
Taking $h=0$ in \eqref{eq:branched-family} and
using $w(e)=1$ yields \eqref{eq:branched-packing}.
\end{proof}

\newpage

It is easy to check that for $(b,c)=(1,0)$, the weighted paths are precisely
those of Lemma~\ref{lem:14325-paths}. In that case,
the additional tridiagonal reduction of
Lemma~\ref{lem:14325-intertwining} gives the
Jacobi fraction of Theorem~\ref{thm:14325}.

\section{`Short' patterns, the original motivating table, and conjectured positivity}
\label{sec:scope}\label{sec:conjclose}

The results of the preceding sections suffice to settle the moment question
for every pattern of length at most four, completing the length-four
classification initiated in \cite{BlitvicKammounSteingrimsson2023}.

\begin{corollary}\label{cor:small}
For $K\le4$, every pattern $\pi\in S_K$ and every admissible overlap
$1\le m<K$ give a Stieltjes moment sequence
$(N_\pi^m(r+1))_{r\ge0}$. The same holds for separated occurrences,
i.e. $m\le0$. The convention $m=K$ corresponds to the law $\delta_1$.
\end{corollary}

\begin{proof}
Proposition~\ref{prop:disjoint} covers $m\le0$ and $m=K$.
At overlap one, Proposition~\ref{prop:extreme-endpoint} and
Theorem~\ref{thm:moments} apply, since, for $K\le4$, $\pi_1$ and 
$\pi_K$ are either adjacent or include $1$ or $K$.
If $m=K-1$, admissibility forces all consecutive comparisons
in $\pi$ to agree, so the pattern is monotone and its packing
posets are chains.

The only remaining case is $(K,m)=(4,2)$. Up to the equivalences in
Part~I~\cite[Section~3]{BKS1} the only relevant prefixes are $12,14$, both of which give chains, and $13$, which gives the Catalan sequence covered by Corollary~\ref{cor:fuss}.
\end{proof}

The preceding results also account for nearly all
the entries in the enumeration table
\cite[Table~1]{BlitvicKammounSteingrimsson2023}
that motivated this work. Table~\ref{tab:moment-status}
records the corresponding positivity results for all the patterns
of length up to five at admissible overlaps $1\le m<K$,
after accounting for the equivalences in
Part~I~\cite[Section~3]{BKS1}. There are fourteen distinct
sequences, thirteen of which are now solved.

For completeness, at overlap three, we have the following known cases.

\begin{proposition}\label{prop:short-overlap-three}
At overlap $m=3$, the patterns $13254$, $14253$, and $15243$
have packing counts
\[
N_{13254}^3(j)=N_{15243}^3(j)=1,
\qquad
N_{14253}^3(j)=C_j
\qquad(j\ge1).
\]
Their moment sequences are represented by $\delta_1$ for
$13254$ and $15243$, and by $x\mu_1(dx)$ for $14253$,
where $\mu_1$ is the $p=1$ Fuss--Catalan distribution
defined in Section~\ref{sec:grids}.
\end{proposition}

\begin{proof}
The two constant enumerations were already recorded in
\cite[Table~1]{BlitvicKammounSteingrimsson2023}.
The Catalan enumeration follows from
\cite[Proposition~2]{BlitvicKammounSteingrimsson2023},
since $31425$ is the reverse complement of $14253$.
The indexing in that proposition counts occurrences beyond
the first, so its $C_{j+1}$ becomes $C_j$ here.

The constant sequences have representing measure $\delta_1$.
For $14253$, the $r$th moment of $x\mu_1(dx)$ is
$C_{r+1}$ by Corollary~\ref{cor:fuss}.
\end{proof}

Furthermore, by the prefix--suffix exchange in
Part~I~\cite[Section~3]{BKS1}, the patterns $13524$ and
$24513$ at overlap two give the same remaining sequence, namely:
\[
\begin{aligned}
\gamma_n:=N_{13524}^{2}(n+1)=N_{24513}^{2}(n+1)
={}&1,\,6,\,81,\,1806,\,57447,\,2399058,\\
&124771320,\,7794139680,\ldots\qquad(n\ge0).
\end{aligned}
\]
For this sequence, our general recurrence in
\cite[Theorem~2]{BlitvicKammounSteingrimsson2023}
reduces to the following triangular form.

\begin{proposition}\label{prop:gamma-recurrence}
Define a triangular array by $a_{0,0}=1$ and
\begin{equation}\label{eq:gamma-recurrence}
a_{n,k}
=
a_{n,k-1}+(2n+k+1)a_{n-1,k}
\qquad(n\ge1,\ 0\le k\le n),
\end{equation}
where $a_{n,k}=0$ when $k<0$ or $k>n$. Then
\begin{equation}\label{eq:gamma-row-sum}
\gamma_n
=\sum_{k=0}^{n}a_{n,k}
\qquad(n\ge0).
\end{equation}
\end{proposition}

\begin{proof}
For $24513$, the suffix recurrence starts from the terminal
overlap $(1,3)$. When another occurrence is appended, an old
overlap $(1,r)$ becomes $(2,r+2)$, and the new terminal overlap
must be $(1,b)$ with $2<b<r+2$. Inductively, every terminal
overlap therefore has the form $(1,b)$ with $3\le b\le n+3$
for a packing with $n+1$ occurrences.

Let $g_n(b)$ count these packings with terminal overlap $(1,b)$,
and set $g_n(b)=0$ outside this range. The new middle letter
must exceed $r+2$, giving $3n+3-r$ choices in a permutation
of length $3n+5$. Thus equation~(2.1) of the cited theorem
specializes to
\[
g_n(b)
=
\sum_{r=b-1}^{n+2}(3n+3-r)g_{n-1}(r)
\qquad(n\ge1,\ 3\le b\le n+3),
\]
with $g_0(3)=1$.

Taking successive differences gives
\[
g_n(n+3-k)-g_n(n+4-k)
=
(2n+k+1)g_{n-1}(n+2-k)
\qquad(0\le k\le n).
\]
Hence $g_n(n+3-k)$ satisfies \eqref{eq:gamma-recurrence}
with the same initial and boundary values as $a_{n,k}$,
so
\[
a_{n,k}=g_n(n+3-k).
\]
Summing over the terminal overlap gives
\[
\gamma_n
=
N_{24513}^{2}(n+1)
=
\sum_{b=3}^{n+3}g_n(b)
=
\sum_{k=0}^{n}a_{n,k},
\]
as required.
\end{proof}

\begin{remark}
Comparison of \eqref{eq:gamma-recurrence} with unweighted ballot
paths gives $\gamma_n\le 3^n n!C_{n+1}$.
\end{remark}

Computations using the first $2000$ terms found both
$H_k(\gamma)$ and $H_k^{(1)}(\gamma)$ to be positive
semidefinite for every tested $k$.
It is therefore reasonable to assume that packings of any pattern of
length up to five give a Stieltjes moment sequence at every
admissible overlap.

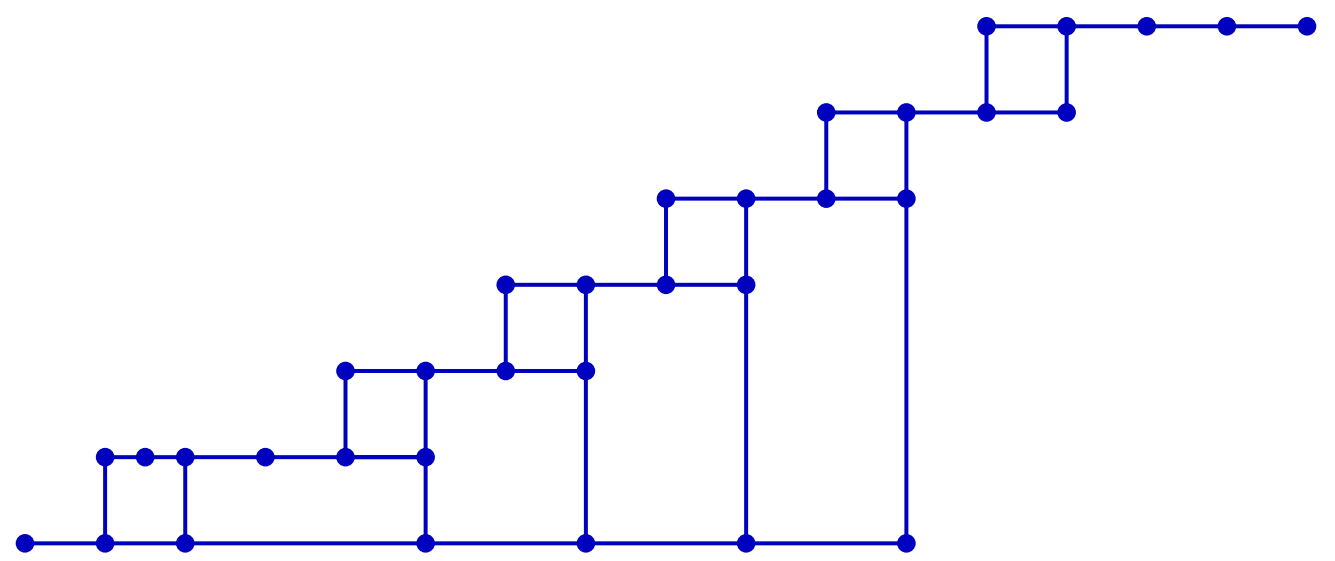
\begin{figure}[htbp]
\centering
\begin{tikzpicture}[x=5.3mm,y=5.7mm,
  edge/.style={draw=blue!75!black,line width=.75pt},
  vertex/.style={circle,fill=blue!75!black,draw=none,
                 inner sep=0pt,minimum size=3.5pt}]
  \draw[edge] (-1,0)--(10,0);
  \foreach \x in {-1,0,1,4,6,8,10}
    \node[vertex] at (\x,0) {};
  \draw[edge] (0,0)--(0,1)--(4,1);
  \draw[edge] (1,0)--(1,1);
  \foreach \x in {0,.5,1,2}
    \node[vertex] at (\x,1) {};
  \foreach \t in {0,...,4} {
    \pgfmathtruncatemacro{\xx}{3+2*\t}
    \pgfmathtruncatemacro{\yy}{1+\t}
    \draw[edge] (\xx,\yy) rectangle
                ({\xx+1},{\yy+1});
    \foreach \dx in {0,1}
      \foreach \dy in {0,1}
        \node[vertex] at ({\xx+\dx},{\yy+\dy}) {};
  }
  \foreach \t in {0,...,3} {
    \pgfmathtruncatemacro{\xx}{4+2*\t}
    \pgfmathtruncatemacro{\yy}{2+\t}
    \draw[edge] (\xx,\yy)--({\xx+1},\yy);
  }
  \foreach \x/\y in {4/1,6/2,8/3,10/4}
    \draw[edge] (\x,0)--(\x,\y);
  \draw[edge] (12,6)--(15,6);
  \foreach \x in {13,14,15}
    \node[vertex] at (\x,6) {};
\end{tikzpicture}
\caption{The poset for six occurrences of $134682579$
at overlap four. Edges are directed to the right or upwards.}
\label{fig-134682579}
\end{figure}

The general conjecture, formulated in
\cite[Conjecture~1]{BlitvicKammounSteingrimsson2023},
according to which $(N_\pi^m(r+1))_{r\ge0}$ is a Hamburger
moment sequence for every pattern $\pi$ and every admissible
overlap $m$, is easily invalidated. For example, it suffices
to look at the pattern $\pi=134682579$ at overlap four
(see Figure~\ref{fig-134682579}), whose leading $5\times5$
Hankel determinant is negative:
\[
\det\bigl(N_{134682579}^{4}(i+j+1)\bigr)_{i,j=0}^{4}
=-63986350161920<0.
\]

Yet, many sequences appear to become Hamburger after discarding
a few initial terms. In the preceding example, all tested
Hankel determinants beginning at occurrence two or later are
positive within the first $3{,}000$ terms.

From our general results, we know that
Theorem~\ref{thm:rooted-tails} guarantees a Stieltjes tail
for every rooted-tree family with
$1\le m\le K/2$, beginning no later than
occurrence $\max(1,\lfloor m/2\rfloor)$. In particular,
the original sequence is already Stieltjes when $m\le3$.

For grids, Theorem~\ref{thm:grid} proves Stieltjes positivity
from the first term for every two-level family of shifted grids,
with arbitrary horizontal and vertical subdivision parameters.
Theorem~\ref{thm:135246} also gives a Stieltjes law for the
three-level grid associated with $135246$ at overlap three.
Positivity for arbitrary numbers of levels remains conjectural.

Furthermore, a more general picture linking poset structure and positivity remains murky. Indeed, with the help of LLMs, we have falsified
several budding structural conjectures linking positivity
to notions of poset `regularity'. And yet general results
are manifestly possible.

Following extensive numerical verification, we are allowing
ourselves the following conjectures.

\begin{conjecture}[Shifted grids]\label{conj:canonical-grids}
Let $\pi\in S_K$ have an admissible overlap
$2\le m\le K/2$, and let
$p_1<\cdots<p_m$ and $s_1<\cdots<s_m$ be the sorted
letters in its $m$-prefix and $m$-suffix, respectively.
If
\[
1=p_1<s_1<p_2<s_2<\cdots<p_m<s_m=K,
\]
then $(N_\pi^m(r+1))_{r\ge0}$ is a Stieltjes moment sequence.
\end{conjecture}

\begin{conjecture}[Overlap one]\label{conj:overlap-one}
For every $K\ge2$ and every $\pi\in S_K$, the sequence
$(N_\pi^1(r+1))_{r\ge0}$ is a Stieltjes moment sequence.
\end{conjecture}

\section*{Acknowledgements}

We thank Cyril Banderier for pointing out the connection between our probabilistic approach to enumeration and the density method.

During the final stages of this project, OpenAI’s GPT-6 assisted with the development of Section~\ref{sec:rooted}, parts of the proofs
of Theorems~\ref{thm:grid}, \ref{thm:gap-two}, and~\ref{thm:14325}, and 
an ansatz that simplified the proof of Theorem~\ref{thm:135246}. It further assisted with critical review and refinement of the manuscript's exposition. 

The research questions as well as the vast majority of the work, including the underlying combinatorial and probabilistic frameworks, were developed by the authors and in the prior (cited) works, and presented to the community over a number of years.

\bibliographystyle{amsplain}
\bibliography{BKS2_authorial_revision_full}

@misc{brignall2026,
      title={Shuffle automata and the growth of 1324-avoiding permutations}, 
      author={Robert Brignall},
      year={2026},
      eprint={2610.02535},
      archivePrefix={arXiv},
      primaryClass={math.CO},
      url={https://arxiv.org/abs/2610.02535}, 
}

@misc{BKS1,
      title={Permutations, Pattern Covers, Posets and Positivity: Part I: Patterns and Posets}, 
      author={Natasha Blitvić and Mohammed Slim Kammoun and Einar Steingrímsson},
      year={2026},
      Note={arXiv:2610.29846},
      archivePrefix={arXiv},
      primaryClass={math.CO},
      url={https://arxiv.org/abs/2610.29846}
}

@article{BlitvicSteingrimsson2021,
  author = {Blitvi{\'c}, Natasha and Steingr{\'i}msson, Einar},
  title = {Permutations, moments, measures},
  journal = {Transactions of the American Mathematical Society},
  volume = {374},
  number = {8},
  pages = {5473--5508},
  year = {2021},
  doi = {10.1090/tran/8330}
}

@article{BanderierWallner2021,
  author = {Banderier, Cyril and Wallner, Michael},
  title = {Young tableaux with periodic walls: counting with the density method},
  journal = {S{\'e}minaire Lotharingien de Combinatoire},
  volume = {85B},
  pages = {Article 47},
  note = {12 pp.},
  year = {2021}
}

@inproceedings{AlbenqueBouttier2012,
  author    = {Albenque, Marie and Bouttier, J{\'e}r{\'e}mie},
  title     = {Constellations and multicontinued fractions:
               Application to {Eulerian} triangulations},
  booktitle = {24th International Conference on Formal Power
               Series and Algebraic Combinatorics (FPSAC 2012)},
  series    = {Discrete Mathematics \& Theoretical Computer
               Science Proceedings},
  volume    = {AR},
  pages     = {805--816},
  year      = {2012},
  doi       = {10.46298/dmtcs.3084}
}

@article{PetreolleSokalZhu2023,
  author = {P{\'e}tr{\'e}olle, Mathias and
            Sokal, Alan D. and
            Zhu, Bao-Xuan},
  title = {Lattice paths and branched continued fractions:
           An infinite sequence of generalizations of the
           Stieltjes--Rogers and Thron--Rogers polynomials,
           with coefficientwise Hankel-total positivity},
  journal = {Memoirs of the American Mathematical Society},
  volume = {291},
  number = {1450},
  year = {2023},
  doi = {10.1090/memo/1450}
}

@book{ST43,
  author    = {J. A. Shohat and J. D. Tamarkin},
  title     = {The Problem of Moments},
  series    = {Mathematical Surveys},
  volume    = {1},
  publisher = {American Mathematical Society},
  year      = {1943}
}

@article{SokalZeng2022,
  author  = {Sokal, Alan D. and Zeng, Jiang},
  title   = {Some multivariate master polynomials for permutations,
             set partitions, and perfect matchings,
             and their continued fractions},
  journal = {Adv. in Appl. Math.},
  volume  = {138},
  year    = {2022},
  pages   = {102341},
  doi     = {10.1016/j.aam.2022.102341}
}

@article{ElizaldeNoy2012,
  author = {Elizalde, Sergi and Noy, Marc},
  title = {Clusters, generating functions and asymptotics for consecutive patterns in permutations},
  journal = {Advances in Applied Mathematics},
  volume = {49},
  number = {3--5},
  pages = {351--374},
  year = {2012}
}

@article{BostanEtAl2020,
 author={Bostan, Alin and Elvey Price, Andrew and Guttmann, Anthony John and Maillard, Jean-Marie},
 title={{Stieltjes} moment sequences for pattern-avoiding permutations},
 journal={The Electronic Journal of Combinatorics}, volume={27}, number={4}, year={2020},
 pages={P4.20}, doi={10.37236/9402}}

@article{BlitvicKammounSteingrimsson2023,
 author={Blitvi{\'c}, Natasha and Kammoun, Mohammed Slim and Steingr{\'i}msson, Einar},
 title={A new perspective on positivity in (consecutive) permutation patterns},
 journal={S{\'e}minaire Lotharingien de Combinatoire},
 volume={89B}, year={2023}, pages={Article 61}, note={12 pp.}}

@article{DwyerElizalde2018,
 author={Dwyer, Tim and Elizalde, Sergi},
 title={{Wilf} equivalence relations for consecutive patterns},
 journal={Advances in Applied Mathematics}, volume={99}, year={2018}, pages={134--157}}

@article{Sokal2023,
 author={Sokal, Alan D.},
 title={Total positivity of some polynomial matrices that enumerate labeled trees and forests {I}: forests of rooted labeled trees},
 journal={Monatshefte f{\"u}r Mathematik}, volume={200}, year={2023}, pages={389--452},
 doi={10.1007/s00605-022-01687-0}}

@article{Sun2017,
  author  = {Sun, Ping},
  title   = {Enumeration of standard Young tableaux
             of shifted strips with constant width},
  journal = {Electronic Journal of Combinatorics},
  volume  = {24},
  number  = {2},
  pages   = {P2.41},
  year    = {2017},
  doi     = {10.37236/6466}
}

@article{Flajolet1980,
 author={Flajolet, Philippe}, title={Combinatorial aspects of continued fractions},
 journal={Discrete Mathematics}, volume={32}, year={1980}, pages={125--161}}

@article{PensonZyczkowski2011,
 author={Penson, Karol A. and {\.{Z}}yczkowski, Karol},
 title={Product of {Ginibre} matrices: {Fuss--Catalan} and {Raney} distributions},
 journal={Physical Review E}, volume={83}, year={2011}, pages={061118},
 doi={10.1103/PhysRevE.83.061118}}
\end{document}